\documentclass[oneside]{scrartcl} 
\usepackage{amsmath,amsthm,amscd,amsfonts}
\numberwithin{equation}{section}
\theoremstyle{definition}
\newtheorem{dfn}{Definition}[section]
\newtheorem{example}[dfn]{Example}
\newtheorem{rem}[dfn]{Remark}
\theoremstyle{plain}
\newtheorem{thm}[dfn]{Theorem}
\newtheorem{prop}[dfn]{Proposition}
\newtheorem{cor}[dfn]{Corollary}
\newtheorem{lem}[dfn]{Lemma}

\title{The elliptic $GL(n)$ dynamical quantum group as an $\mathfrak{h}$-Hopf algebroid}
\author{Jonas T. Hartwig\footnote{The author was partially supported by the
Netherlands Organization for Scientific Research (NWO) in the
VIDI-project ``Symmetry and modularity in exactly solvable models''.}}
\date{ }

\def\al{\alpha}
\def\be{\beta}
\def\ga{\gamma}
\def\de{\delta}\def\De{\Delta}
\def\ep{\varepsilon}
\def\ze{\zeta}

\def\th{\theta}

\def\la{\lambda} \def\La{\Lambda}
\def\rh{\rho}
\def\si{\sigma}
\def\ta{\tau}
\def\ph{\varphi}

\def\om{\omega}
\def\hf{\mathfrak{h}}
\def\C{\mathbb{C}} \def\Z{\mathbb{Z}} 
\def\R{\mathsf{R}} \def\L{\mathsf{L}}
 \def\Cdd{W} \def\Che{W} \def\Cher{\mathsf{C}}
\def\Fell{\mathcal{F}_{\mathrm{ell}}} \def\Fscr{\mathcal{F}}
\def\Id{\mathrm{Id}}
\def\oneton{[1,n]} \def\onetod{[1,d]}

\newcommand{\sdeg}{\mathrm{sdeg}}
\DeclareMathOperator{\tr}{tr} 
\DeclareMathOperator{\diag}{diag} 
\DeclareMathOperator{\im}{Im}
\DeclareMathOperator{\End}{End}
\DeclareMathOperator{\sgn}{sgn}

\begin{document}
\maketitle
\begin{abstract}
Using the language of $\hf$-Hopf algebroids which was
introduced by Etingof and Varchenko,
we construct a dynamical quantum group, $\Fell(GL(n))$,
from the elliptic 
solution of the quantum dynamical Yang-Baxter
equation with spectral parameter associated to the
Lie algebra $\mathfrak{sl}_n$.
We apply the generalized FRST construction
and obtain an $\hf$-bialgebroid $\Fell(M(n))$.
Natural analogs of
the exterior algebra and their matrix elements, elliptic minors,
are defined and studied.
We show how to use the cobraiding to prove that the elliptic determinant
is central.
Localizing at this determinant and constructing an antipode we obtain
the $\hf$-Hopf algebroid $\Fell(GL(n))$. 
\end{abstract}
\tableofcontents

\section{Introduction}

The quantum dynamical Yang-Baxter (QDYB) equation
was introduced by Gervais and Neveu \cite{GeNe}.
It was realized by Felder \cite{Fe95a} that this equation is equivalent to
the Star-Triangle relation in statistical mechanics.
It is a generalization of the quantum Yang-Baxter equation,
involving an extra, so called dynamical, parameter. 
In \cite{Fe95a} an interesting elliptic solution to the QDYB equation with
spectral parameter was given, adapted from
the $A_n^{(1)}$ solution to the Star-Triangle relation constructed in \cite{JKMO88}.
Felder also defined a tensor category,
which he suggested should be thought of as an elliptic analog
of the category of representations of quantum groups.
This category was further studied in \cite{FeVa96} in
the $\mathfrak{sl}_2$ case.

In \cite{FeVa97}, the authors considered objects in Felder's category which
were proposed as analogs of exterior and symmetric powers of the vector
representation of $\mathfrak{gl}_n$.
To each object in the tensor category they associate an algebra of
vector-valued difference operators and prove that a certain operator,
constructed from the analog of the top exterior power, commutes with all
other difference operators.
This is also proved in \cite{TaVa01} (Appendix B) in more detail and
in \cite{ZhShYu} using a different approach.

An algebraic framework for studying dynamical R-matrices without
spectral parameter was introduced
in \cite{EtVa98}. There the authors defined the notion of $\hf$-bialgebroids
and $\hf$-Hopf algebroids, a special case of the Hopf algebroids defined
by Lu \cite{Lu96}. They also show, using a generalized version of
the FRST construction, how to associate to every solution $R$ of the
non-spectral quantum dynamical Yang-Baxter equation an $\hf$-bialgebroid.
Under some extra condition
they get an $\hf$-Hopf algebroid by adjoining
formally the matrix elements of the inverse L-matrix.
This correspondence gives a tensor equivalence between the category of representations
of the R-matrix and the category of so called dynamical representations
of the $\hf$-bialgebroid.

In this paper we define an $\hf$-Hopf algebroid associated to
the elliptic R-matrix from \cite{Fe95a}
with both dynamical and spectral parameter for $\mathfrak{g}=\mathfrak{sl}_n$.
This generalizes the
spectral elliptic dynamical $GL(2)$ quantum group from \cite{KovNRo}
and the non-spectral trigonometric dynamical $GL(n)$ quantum group
from \cite{KovN}.  
As in \cite{KovNRo}, this is done by first using the
the generalized FRST construction, modified to also include spectral parameters.
In addition to the usual
RLL-relation, residual relations must be added ``by hand''
to be able to prove that different expressions for the determinant are equal.

Instead of
adjoining formally all the matrix elements of the inverse L-matrix,
we adjoin only the inverse of the determinant, as in \cite{KovNRo}.
Then we express the antipode using this inverse.
The main problem is to find the correct formula for the determinant,
to prove that it is central, and to provide row and column expansion
formulas for the determinant in the setting of $\hf$-bialgebroids.

The plan of this paper is as follows.
After introducing some notation in Section \ref{sec:notation},
we recall the definition of the elliptic R-matrix in
Section \ref{sec:theRmatrix}.
In Section \ref{sec:bialgebroids} we review the definition of
$\hf$-bialgebroids and
the generalized FRST construction
with special emphasize on how to treat residual relations for
a general R-matrix.
We write down the relations explicitly in Section \ref{sec:thealgFellMn}
for the algebra $\Fell(M(n))$ obtained
from the elliptic R-matrix. In particular we show
that only one family of residual identities are needed.

Left and right analogs of the exterior algebra over $\C^n$ is defined in
Section \ref{sec:extalgs}
in a similar way as in \cite{KovN}. They are certain comodule algebras over
$\Fell(M(n))$ and arise naturally from a single relation
analogous to $v\wedge v=0$.
The matrix elements of these corepresentations
are generalized minors depending on a spectral parameter.
Their properties are studied
in Section \ref{sec:minors}.
In particular we show that the left and right versions
of the minors in fact coincide.
In Section \ref{sec:laplace} we prove Laplace expansion formulas for
these elliptic quantum minors.

In Section \ref{sec:pairingsanddet} we show that the $\hf$-bialgebroid
$\Fell(M(n))$ can be equipped with a
cobraiding, in the sense of \cite{Ro}. We use this and
the the ideas as in \cite{FeVa97} and \cite{TaVa01}
to prove that the determinant is central for all values of the spectral
parameters.
This implies that the determinant is central in the operator algebra
as shown in \cite{FeVa97}.

Finally, in Section \ref{sec:antipode} we define
$\Fell(GL(n))$ to be the localization of $\Fell(M(n))$ at the determinant
and show that it has an antipode giving it the structure
of an $\hf$-Hopf algebroid.

\section{Preliminaries}
\subsection{Notation} \label{sec:notation}
Let $p,q\in\mathbb{R}$, $0<p,q<1$.
We assume $p,q$ are generic in the sense that
if $p^aq^b=1$ for some $a,b\in\Z$, then $a=b=0$.

Denote by $\th$ the normalized Jacobi theta function:
\begin{equation}\label{eq:thetadef}
\theta(z)=\theta(z;p)=\prod_{j=0}^\infty (1-zp^j)(1-p^{j+1}/z).
\end{equation}
It is holomorphic on $\C^\times:=\C\backslash\{0\}$ with zero-set $\{p^k:k\in\Z\}$ and satisfies
\begin{equation}\label{eq:thetaid}
\th(z^{-1})=\th(pz)=-z^{-1}\th(z)
\end{equation}
and the addition formula
\begin{equation}\label{eq:thetaadd}
\th(xy,x/y,zw,z/w)=\th(xw,x/w,zy,z/y)+(z/y)\th(xz,x/z,yw,y/w),
\end{equation}
where we use the notation
\[\theta(z_1,\ldots,z_n)=\theta(z_1)\cdots\theta(z_n).\]
Recall also the Jacobi triple product identity, which can be written
\begin{equation}
  \label{eq:jacobitriple}
  \sum_{k\in\Z} (-z)^{k} p^{\frac{k(k-1)}{2}} = \th(z)\prod_{j=1}^\infty (1-p^j).
\end{equation}
It will sometimes be convenient to use the auxiliary function $E$ given by
\begin{equation} \label{eq:Edef}
 E:\C\to\C,\qquad E(s)=q^s\th(q^{-2s}).
\end{equation}
Relation \eqref{eq:thetaid} implies that $E(-s)=-E(s)$.

The set $\{1,2,\ldots,n\}$ will be denoted by $\oneton$.

\subsection{The elliptic R-matrix} \label{sec:theRmatrix}
Let $\hf$ be a complex vector space, viewed as an abelian Lie algebra,
$\hf^*$ its dual space
and $V=\bigoplus_{\la\in\hf^*} V_\la$ a diagonalizable $\hf$-module.
A dynamical R-matrix is by definition a meromorphic function
\[R:\hf^*\times \C^\times \to \End_\hf(V\otimes V)\]
satisfying the quantum dynamical Yang-Baxter equation with spectral parameter (QDYBE):
\begin{multline} \label{eq:qdybe}
R(\la,\frac{z_2}{z_3})^{(23)}R(\la-h^2,\frac{z_1}{z_3})^{(13)}R(\la,\frac{z_1}{z_2})^{(12)}=\\
=R(\la-h^3,\frac{z_1}{z_2})^{(12)}R(\la,\frac{z_1}{z_3})^{(13)}R(\la-h^1,\frac{z_2}{z_3})^{(23)}.
\end{multline}
Equation \eqref{eq:qdybe} is an equality in the algebra of meromorphic functions
$\hf^*\times\C^\times\to \End(V^{\otimes 3})$. Upper indices are
leg-numbering notation and $h$ indicates the action of $\hf$. For example,
\[R(\la-h^3,\frac{z_1}{z_2})^{(12)}(u\otimes v\otimes w) = R(\la-\al,\frac{z_1}{z_2})(u\otimes v)
\otimes w,\qquad\text{if $w\in V_\al$.}\]
An R-matrix $R$ is called \emph{unitary} if
\begin{equation}\label{eq:Runitary}
R(\la,z)R(\la,z^{-1})^{(21)}=\Id_{V\otimes V}
\end{equation}
as meromorphic functions on $\hf^*\times \C^{\times}$ with
values in $\End_{\hf}(V\otimes V)$.

In the example we study, $\hf$ is 
the Cartan subalgebra of $\mathfrak{sl}(n)$. Thus $\hf$ is
the abelian Lie algebra of all
traceless diagonal complex $n\times n$ matrices.
Let $V$ be the $\hf$-module $\C^n$ with standard basis $e_1,\ldots,e_n$.
Define $\om(i)\in\hf^*$ ($i=1,\ldots,n$) by 
\[\om(i)(h)=h_i,\quad \text{if }h=\diag(h_1,\ldots,h_n)\in\hf.\]
We have $V=\bigoplus_{i=1}^n V_{\om(i)}$ and
$V_{\om(i)}=\C e_i$.
Define \[R:\hf^*\times\C^\times\to\End(V\otimes V)\] by
\begin{equation}
  \label{eq:Rmatrixdef}
R(\la,z)=\sum_{i=1}^n E_{ii}\otimes E_{ii}+
\sum_{i\neq j}\al(\la_{ij},z)E_{ii}\otimes E_{jj}+
\sum_{i\neq j}\be(\la_{ij},z)E_{ij}\otimes E_{ji}
\end{equation}
where $E_{ij}\in\End(V)$ are the matrix units, $\la_{ij}$ ($\la\in\hf^*$) is an abbrevation for
$\la(E_{ii}-E_{jj})$, and
\[\al,\be :\C\times\C^\times\to \C\]
are given by
\begin{equation}\label{eq:aldef}
\alpha(\la,z)=\al(\la,z;p,q)=\frac{\th(z)\th(q^{2(\la+1)})}{\th(q^2z)\th(q^{2\la})},
\end{equation}
\begin{equation}\label{eq:bedef}
\beta(\la,z)=\be(\la,z;p,q)=\frac{\th(q^2)\th(q^{-2\la}z)}{\th (q^2z)\th (q^{-2\la})}.
\end{equation}

\begin{prop}[\cite{Fe95a}]\label{prop:RisanRmatrix}
The map $R$ is a unitary R-matrix.
\end{prop}

For the readers convenience, we give
the explicit relationship between the R-matrix \eqref{eq:Rmatrixdef}
and Felders R-matrix as written in \cite{FeVa97} which we denote by $R_1$.
Thus $R_1:\hf_1^*\times\C\to\End(V\otimes V)$,
where $\hf_1$ is the Cartan subalgebra of $\mathfrak{gl}(n)$,
is defined as in \eqref{eq:Rmatrixdef} with $\al,\be$
replaced by $\al_1,\be_1:\C^2\to\C$,
\begin{align}
\label{eq:alpha1def}
\al_1(\la,x)&=
\al_1(\la,x;\ta,\ga)=\frac{\th_1(x;\ta)\th_1(\la+\ga;\ta)}{\th_1(x-\ga;\ta)\th_1(\la;\ta)},\\
\label{eq:beta1def}
\be_1(\la,x)&=
\be_1(\la,x;\ta,\ga)=-\frac{\th_1(x+\la;\ta)\th_1(\ga;\ta)}{\th_1(x-\ga;\ta)\th_1(\la;\ta)}.
\end{align}
Here $\ta,\ga\in\C$ with $\im\ta>0$ and $\th_1$ is the first Jacobi theta function
\[\th_1(x;\ta)=-\sum_{j\in \Z+\frac{1}{2}} e^{\pi ij^2 \ta+2\pi ij(x+1/2)}.\]
As proved in \cite{Fe95a}, $R_1$ satisfies the following version of the QDYBE:
\begin{multline}
  \label{eq:qdybediff}
  R_1(\la-\ga h^3,x_1-x_2)^{(12)}R_1(\la,x_1-x_3)^{(13)}R_1(\la-\ga h^1,x_2-x_3)^{(23)}=\\=
R_1(\la,x_2-x_3)^{(23)}R_1(\la-\ga h^2,x_1-x_3)^{(13)}R_1(\la,x_1-x_2)^{(12)}
\end{multline}
and the unitarity condition
\begin{equation}
  \label{eq:unitaritydiff}
  R_1(\la,x)R_1^{21}(\la,-x)=\Id_{V\otimes V}.
\end{equation}
We can identify $\hf^*\simeq \hf_1^*/\C\tr$ where $\tr\in \hf_1^*$ is the trace.
Since $R_1$ has the form \eqref{eq:Rmatrixdef}, it is constant, as a function of $\la\in\hf_1^*$,
on the cosets modulo $\C\tr$. So $R_1$ induces a map $\hf^*\times\C\to\End(V\otimes V)$, which
we also denote by $R_1$, still satisfying \eqref{eq:qdybediff},\eqref{eq:unitaritydiff}.

Let $\ta,\ga\in\C$ with $\im\ta>0$ be such that $p=e^{\pi i\ta}$, $q=e^{\pi i\ga}$.
Then, as meromorphic functions of $(\la,x)\in\hf^*\times\C$,
\begin{equation}
  \label{eq:R1Rrelation}
R_1(\ga \la, -x; \ta/2, \ga)=R(\la,z; p,q)
\end{equation}
where $z=e^{2\pi ix}$.
Indeed, using the Jacobi triple product identity \eqref{eq:jacobitriple} we have
\[\th_1(x;\ta/2)=ie^{\pi i (\ta/2-x)} \th(z)\prod_{j=1}^\infty (1-p^j), \]
and substituting this into \eqref{eq:alpha1def} and \eqref{eq:beta1def} gives
$\al_1(\ga\la,-x;\ta/2,\ga)=\al(\la,z;p,q)$ and
$\be_1(\ga\la,-x;\ta/2,\ga)=\be(\la,z;p,q)$ which proves \eqref{eq:R1Rrelation}.

By replacing $\la$, $x_i$ in \eqref{eq:qdybediff}
by $\ga\la$, $-x_i$ and using \eqref{eq:R1Rrelation} 
we obtain \eqref{eq:qdybe} with $z_i=e^{2\pi ix_i}$. Similarly
the unitarity \eqref{eq:Runitary} of $R$ is
obtained from \eqref{eq:unitaritydiff}.

\subsection{Useful identities}
We end this section by recording some useful identities.
Recall the definitions of $\al,\be$ in \eqref{eq:aldef},\eqref{eq:bedef}.
It is immediate that
\begin{equation}\label{eq:albeq2}
\al(\la,q^2)=\be(-\la,q^2).
\end{equation}
By induction, one generalizes \eqref{eq:thetaid} to
\begin{equation}\label{eq:leftel3}
\th(p^sz)=(-1)^s(p^{s(s-1)/2}z^s)^{-1}\th(z),\quad\text{for }s\in\Z.
\end{equation}
Applying \eqref{eq:leftel3} to the definitions of $\al$,$\be$ we get
\begin{equation}\label{eq:albepk}
\al(\la,p^kz)=q^{2k}\al(\la,z),\quad \be(\la,p^kz)=q^{2k(\la+1)}\be(\la,z),
\end{equation}
and, using also $\th(z^{-1})=-z^{-1}\th(z)$,
\begin{equation}\label{eq:albelimits}
\begin{split}
&\lim_{z\to p^{-k}q^{-2}}\frac{q^{-1}\th(q^2z)}{q\th(q^{-2}z)}\al(\la,z)=\al(\la,p^kq^2),\\
&\lim_{z\to p^{-k}q^{-2}}\frac{q^{-1}\th(q^2z)}{q\th(q^{-2}z)}\be(\la,z)=-\be(-\la,p^kq^2),
\end{split}
\end{equation}
for $\la\in\C$, $z\in\C^\times$, and $k\in\Z$.
By the addition formula \eqref{eq:thetaadd} with
\[(x,y,z,w)=(z^{1/2}q^{-\la+1},z^{1/2}q^{\la-1},z^{1/2}q^{\la+1},z^{1/2}q^{-\la-1})\]
we have
\begin{equation}\label{eq:albeth}
\al(\la,z)\al(-\la,z)-\be(\la,z)\be(-\la,z)=
q^2\frac{\th(q^{-2}z)}{\th(q^2z)}.
\end{equation}

\section{$\hf$-Bialgebroids} \label{sec:bialgebroids}
\subsection{Definitions}
We recall the some definitions from \cite{EtVa98}.
Let $\hf^*$ be a finite-dimensional complex vector space (for example the dual space
of an abelian Lie algebra) and
$M_{\hf^*}$ be the field of meromorphic functions on $\hf^*$.
\begin{dfn}
 An \emph{$\hf$-algebra} is a complex associative algebra $A$ with $1$ which is
 bigraded over $\hf^*$, $A=\bigoplus_{\al,\be\in\hf^*} A_{\al\be}$, and equipped with two
 algebra embeddings $\mu_l,\mu_r: M_{\hf^*}\to A$, called the left and right moment maps,
 such that
\begin{equation}
 \mu_l(f)a=a\mu_l(T_\al f),\qquad \mu_r(f)a=a\mu_r(T_\be f),
 \quad\text{for $a\in A_{\al\be}$, $f\in M_{\hf^*}$},
\end{equation}
where $T_\al$ denotes the automorphism $(T_\al f)(\ze)=f(\ze+\al)$ of $M_{\hf^*}$.
A morphism of $\hf$-algebras is an algebra homomorphism preserving the bigrading and the
moment maps.
\end{dfn}

The \emph{matrix tensor product} $A\widetilde{\otimes} B$ of two $\hf$-algebras $A$, $B$
is the $\hf^*$-bigraded vector space with
$(A\widetilde{\otimes} B)_{\al\be} = \bigoplus_{\ga\in\hf^*} (A_{\al\ga}\otimes_{M_{\hf^*}} B_{\ga\be})$, 
where $\otimes_{M_{\hf^*}}$ denotes tensor product over $\C$ modulo the relations
\begin{equation}
 \mu_r^A(f)a\otimes b = a \otimes \mu_l^B(f) b,\quad\text{for all $a\in A$, $b\in B$, $f\in M_{\hf^*}$}.
\end{equation}
The multiplication $(a\otimes b)(c\otimes d)=ac\otimes bd$ for $a,c\in A$ and $b,d\in B$ and the
moment maps $\mu_l(f)=\mu_l^A(f)\otimes 1$ and $\mu_r(f)=1\otimes \mu_r^B(f)$ make
$A\widetilde{\otimes} B$ into an $\hf$-algebra.
\begin{example}
 Let $D_\hf$ be the algebra of operators on $M_{\hf^*}$ of the form
 $\sum_i f_iT_{\al_i}$ with $f_i\in M_{\hf^*}$ and $\al_i\in \hf^*$. It is an $\hf$-algebra with bigrading
 $fT_{-\al}\in (D_\hf)_{\al\al}$ and both moment maps equal to the natural embedding.
\end{example}

For any $\hf$-algebra $A$, there are canonical isomorphisms
$A\simeq A\widetilde{\otimes} D_\hf \simeq D_\hf \widetilde{\otimes} A$
defined by
\begin{equation}\label{eq:tensorunit}
 x\simeq x\otimes T_{-\be} \simeq T_{-\al}\otimes x,\quad\text{for $x\in A_{\al\be}$.}
\end{equation}

\begin{dfn}
 An \emph{$\hf$-bialgebroid} is an $\hf$-algebra $A$ equipped with two $\hf$-algebra morphisms,
 the comultiplication $\De:A\to A\widetilde{\otimes} A$ and the counit
 $\ep:A\to D_\hf$ such that $(\De\otimes\Id)\circ\De = (\Id\otimes\De)\circ \De$
 and $(\ep\otimes\Id)\circ\De = \Id = (\Id\otimes\ep)\circ\De$, under the identifications
 \eqref{eq:tensorunit}.
\end{dfn}

\subsection{The generalized FRST-construction}\label{sec:FRST}
In \cite{EtVa98} the authors gave a generalized FRST-construction which
attaches an $\hf$-bialgebroid to each
solution of the quantum dynamical Yang-Baxter equation without
spectral parameter. 
It was described in \cite{KovNRo} one way of extending this
to the case when a spectral parameter is also present. However,
when specifying the R-matrix to \eqref{eq:Rmatrixdef} with $n=2$,
they had to impose in addition certain
so called \emph{residual relations} in order to prove
for example that the determinant is central.
Such relations were also required in \cite{FeVa96} in a different
algebraic setting.
In the setting of operator algebras, where
the algebras consist of linear operators on a vector space
depending meromorphically on the spectral variables, as in
\cite{FeVa97}, such relations are consequences of the ordinary
RLL-relations by taking residues.

Another motivation for our procedure is that
$\hf$-bialgebroids associated to gauge equivalent R-matrices 
should be isomorphic. 
In particular one should be allowed to multiply the R-matrix by
any nonzero meromorphic function of the spectral variable
without changing the isomorphism class of the associated algebra 
(for the full definition of gauge equivalent R-matrices
see \cite{EtVa98}).

These considerations suggest
the following procedure for constructing an $\hf$-bialgebroid from
a quantum dynamical R-matrix with spectral parameter.

Let $\hf$ be a finite-dimensional abelian Lie algebra, $V=\bigoplus_{\al\in\hf^*} V_\al$ a finite-dimensional
diagonalizable $\hf$-module and $R:\hf^*\times \C^\times\to\End_\hf(V\otimes V)$ a meromorphic function.
We attach to this data an \mbox{$\hf$-bialgebroid} $A_R$ as follows.
Let $\{e_x\}_{x\in X}$ be a homogeneous basis of $V$, where $X$ is an index set.
The matrix elements $R_{xy}^{ab}: \hf^*\times\C^\times\to \C$ of $R$ are given by
\begin{equation}
  \label{eq:Rmatrixelements}
R(\ze,z)(e_a\otimes e_b)=\sum_{x,y\in X} R_{xy}^{ab}(\ze,z)e_x\otimes e_y.
\end{equation}
They are meromorphic on $\hf^*\times\C^\times$. Define $\om: X\to\hf^*$ by
$e_x\in V_{\om(x)}$. Let $\tilde A_R$ be the complex associative algebra with $1$
generated by $\{L_{xy}(z): x,y\in X, z\in\C^\times\}$
and two copies of $M_{\hf^*}$, whose elements are denoted by
$f(\la)$ and $f(\rh)$, respectively, with defining relations $f(\la)g(\rh)=g(\rh)f(\la)$
for $f,g\in M_{\hf^*}$ and
\begin{equation}\label{eq:fL}
f(\la)L_{xy}(z)=L_{xy}(z)f(\la+\om(x)),\quad f(\rh)L_{xy}(z)=L_{xy}(z)f(\rh+\om(y)),
\end{equation}
for all $x,y\in X$, $z\in\C^\times$ and $f\in M_{\hf^*}$.
The bigrading on $\tilde A_R$ is given by $L_{xy}(z)\in (\tilde A_R)_{\om(x),\om(y)}$
for $x,y\in X$, $z\in\C^\times$ and $f(\la),f(\rh)\in (\tilde A_R)_{00}$ for $f\in M_{\hf^*}$. 
The moment maps are defined by $\mu_l(f)=f(\la)$, $\mu_r(f)=f(\rh)$.
The counit and comultiplication are defined by
\begin{align}
&\ep(L_{ab}(z))=\de_{ab}T_{-\om(a)},\quad \ep(f(\la))=\ep(f(\rh))=fT_0,\\
&\De(L_{ab}(z))=\sum_{x\in X} L_{ax}(z)\otimes L_{xb}(z),\\
&\De(f(\la))=f(\la)\otimes 1,\quad \De(f(\rh))=1\otimes f(\rh).
\end{align}
This makes $\tilde A_R$ into an $\hf$-bialgebroid.

Consider the ideal in $\tilde A_R$ generated by the RLL-relations 
\begin{equation}  \label{eq:RLLmatrixelements_general}
\sum_{x,y\in X} R_{ac}^{xy}(\la,\frac{z_1}{z_2})L_{xb}(z_1)L_{yd}(z_2)=
\sum_{x,y\in X} R_{xy}^{bd}(\rh,\frac{z_1}{z_2})L_{cy}(z_2)L_{ax}(z_1),
\end{equation}
where $a,b,c,d\in X$, and $z_1,z_2\in\C^\times$.
More precisely, to account for possible singularities of $R$,
we let $I_R$ be the ideal in $\tilde A_R$ generated by all relations of the form
\begin{multline}  \label{eq:RLLmatrixelements_residual_general}
\sum_{x,y\in X} \lim_{w\to z_1/z_2}\left( \ph(w)
R_{ac}^{xy}(\la,w)\right) L_{xb}(z_1)L_{yd}(z_2)=\\=
\sum_{x,y\in X} \lim_{w\to z_1/z_2}\left( \ph(w)
R_{xy}^{bd}(\rh,w)\right) L_{cy}(z_2)L_{ax}(z_1),
\end{multline}
where $a,b,c,d\in X$, $z_1,z_2\in\C^\times$ and
$\ph:\C^\times\to\C$ is a meromorphic function such that
the limits exist.
%

We define $A_R$ to be $\tilde A_R/I_R$. 
The bigrading descends to $A_R$ because \eqref{eq:RLLmatrixelements_residual_general}
is homogeneous,
of bidegree \mbox{$\om(a)+\om(c),\om(b)+\om(d)$}, by the \mbox{$\hf$-invariance} of $R$.
One checks that
$\De(I_R)\subseteq \tilde A_R\widetilde{\otimes} I_R+I_R\widetilde{\otimes}\tilde A_R$
and $\ep(I_R)=0$. Thus $A_R$ is an $\hf$-bialgebroid with the induced maps.

\begin{rem}
Objects in Felder's tensor category associated to an R-matrix $R$
are certain meromorphic functions
$L:\hf^\ast\times\C^\times\to\End_\hf(\C^n\otimes W)$ where $W$
is a finite-dimensional $\hf$-module \cite{Fe95a}. After regularizing $L$
with respect to the spectral parameter it will give rise
to a dynamical representation of the $\hf$-bialgebroid $A_R$
in the same way as in the non-spectral case treated in \cite{EtVa98}.
The residual relations incorporated in
\eqref{eq:RLLmatrixelements_residual_general} are crucial for
this fact to be true in the present, spectral, case.

\end{rem}

\subsection{Operator form of the RLL relations}\label{sec:operatorform}
It is well-known that the RLL-relations \eqref{eq:RLLmatrixelements_general}
can be written as a matrix relation. We show how this is done in
the present setting. It will be used later in Section \ref{sec:cherednik}.

Assume $R_{xy}^{ab}(\ze,z)$ are defined, as
meromorphic functions of $\ze\in\hf^*$, for any $z\in\C^\times$.
Define
$\R(\la,z), \R(\rh,z)\in\End(V\otimes V\otimes A_R)$ by
\begin{align*}
 \R(\la,z)(e_a\otimes e_b\otimes u) &= \sum_{x,y\in X} e_x\otimes e_y\otimes R^{ab}_{xy}(\la,z)u,\\
 \R(\rh,z)(e_a\otimes e_b\otimes u) &= \sum_{x,y\in X} e_x\otimes e_y\otimes R^{ab}_{xy}(\rh,z)u,
\end{align*}
for $a,b\in X$, $u\in A_R$.
Note that the $\la$ and $\rh$ in the left hand side are not variables but merely indicates
which moment map is to be used.
For $z\in\C^\times$ we also define $\L(z)\in\End(V\otimes A_R)$
by
\[\L(z)=\sum_{x,y\in X} E_{xy}\otimes L_{xy}(z).\]
Here $E_{xy}$ are the matrix units in $\End(V)$ and $A_R$ acts on itself by left multiplication.
The RLL relation \eqref{eq:RLLmatrixelements_general} is equivalent to
\begin{equation}\label{eq:RLLoperatorform}
 \R(\la,z_1/z_2)\L^1(z_1)\L^2(z_2)=\L^2(z_2)\L^1(z_1)\R(\rh+h^1+h^2,z_1/z_2)
\end{equation}
in $\End(V\otimes V\otimes A_R)$, where
$\L^i(z)=\L(z)^{(i,3)}\in \End(V\otimes V\otimes A_R)$ for $i=1,2$.
This can be seen by acting on $e_b\otimes e_d\otimes 1$ in
both sides of \eqref{eq:RLLoperatorform}, and collecting and equating
terms of the form $e_a\otimes e_c\otimes u$. The matrix elements of the R-matrix
in the right hand side can then be moved to the left using that $R$ is $\hf$-invariant,
and relation \eqref{eq:fL}.

\section{The algebra $\Fell(M(n))$} \label{sec:thealgFellMn}
We now specialize to the case where
$\hf$ is the Cartan subalgebra of $\mathfrak{sl}(n)$,
$V=\C^n$ and $R$ is given by \eqref{eq:Rmatrixdef}-\eqref{eq:bedef}.
The case $n=2$ was considered in \cite{KovNRo}. 
We will show that \eqref{eq:RLLmatrixelements_residual_general}
contains precisely one additional family of relations, as compared
to \eqref{eq:RLLmatrixelements_general}, and we write down all
relations explicitly.

When we apply the generalized FRST-construction to this data
we obtain an $\hf$-bialgebroid which we denote by $\Fell(M(n))$.
The generators $L_{ij}(z)$ will be denoted by $e_{ij}(z)$.
Thus
$\Fell(M(n))$ is the unital associative $\C$-algebra generated by
$e_{ij}(z)$, $i,j\in\oneton$, $z\in\C^\times$, and two copies of $M_{\hf^*}$,
whose elements are denoted by $f(\la)$ and $f(\rh)$ for $f\in M_{\hf^*}$, subject
to the following relations
\begin{align}\label{eq:FellMnscalarrel}
  f(\la)e_{ij}(z)&=e_{ij}(z)f(\la+\om(i)),&
  f(\rh)e_{ij}(z)&=e_{ij}(z)f(\rh+\om(j)),
\end{align}
for all $f\in M_{\hf^*}$, $i,j\in\oneton$ and $z\in\C^\times$, and
\begin{equation}
  \label{eq:RLLmatrixelements}
\sum_{x,y=1}^n R_{ac}^{xy}(\la,\frac{z_1}{z_2})e_{xb}(z_1)e_{yd}(z_2)=
\sum_{x,y=1}^n R_{xy}^{bd}(\rh,\frac{z_1}{z_2})e_{cy}(z_2)e_{ax}(z_1),
\end{equation}
for all $a,b,c,d\in\oneton$.
More explicitly, from \eqref{eq:Rmatrixdef} we have
\begin{equation}\label{eq:Rmatrixelts}
R^{ab}_{xy}(\ze,z)=\begin{cases}
                    1,&a=b=x=y,\\
		    \al(\ze_{xy},z),&a\neq b, x=a,y=b,\\
		    \be(\ze_{xy},z),&a\neq b, x=b,y=a,\\
		    0,&\text{otherwise,}
                   \end{cases}
\end{equation}
which substituted into \eqref{eq:RLLmatrixelements} yields four families of relations:
\begin{subequations}\label{eq:rel1-4}
\begin{gather}
  \label{eq:rel1}
e_{ab}(z_1)e_{ab}(z_2)=e_{ab}(z_2)e_{ab}(z_1),
\\
  \label{eq:rel2}
  e_{ab}(z_1)e_{ad}(z_2)=\al(\rh_{bd},\frac{z_1}{z_2})e_{ad}(z_2)e_{ab}(z_1)+
  \be(\rh_{db},\frac{z_1}{z_2})e_{ab}(z_2)e_{ad}(z_1),
\\
  \label{eq:rel3}
  \al(\la_{ac},\frac{z_1}{z_2})e_{ab}(z_1)e_{cb}(z_2)+
  \be(\la_{ac},\frac{z_1}{z_2})e_{cb}(z_1)e_{ab}(z_2)=e_{cb}(z_2)e_{ab}(z_1),
\\
\label{eq:rel4}
\begin{split}
\al(\la_{ac},\frac{z_1}{z_2})e_{ab}(&z_1)e_{cd}(z_2)+
\be(\la_{ac},\frac{z_1}{z_2})e_{cb}(z_1)e_{ad}(z_2)=\\
&=\al(\rh_{bd},\frac{z_1}{z_2})e_{cd}(z_2)e_{ab}(z_1)+
\be(\rh_{db},\frac{z_1}{z_2})e_{cb}(z_2)e_{ad}(z_1),
\end{split}
\end{gather}
\end{subequations}
where $a,b,c,d\in\oneton$, $a\neq c$ and $b\neq d$.
Since $\th$ has zeros precisely at $p^k, k\in\Z$, $\al$ and $\be$ have poles
at $z=q^{-2}p^k, k\in\Z$.
Thus \eqref{eq:rel2}-\eqref{eq:rel4}
are to hold for $z_1,z_2\in\C^\times$ with $z_1/z_2\notin\{p^kq^{-2}:k\in\Z\}$.

In \eqref{eq:RLLmatrixelements_residual_general},
assuming $a\neq c$, $b\neq d$, and taking $z_1=z$, $z_2=p^kq^2z$,
$\ph(w)=\frac{q^{-1}\th(q^2w)}{q\th(q^{-2}w)}$,
and using the limit formulas \eqref{eq:albelimits},
we obtain the relation
\begin{multline}\label{eq:rel5}
\al(\la_{ac},q^2)\big(e_{ab}(z)e_{cd}(p^kq^2z)-q^{2k\la_{ca}}e_{cb}(z)e_{ad}(p^kq^2z)\big)=\\
=\al(\rh_{bd},q^2)e_{cd}(p^kq^2z)e_{ab}(z)-
 q^{2k\rh_{bd}}\be(\rh_{bd},q^2)e_{cb}(p^kq^2z)e_{ad}(z).
\end{multline}
This identity does not follow from \eqref{eq:rel1}-\eqref{eq:rel4} in an obvious way.
It will be called the \emph{residual RLL relation}.

\begin{prop}
Relations \eqref{eq:rel1-4},\eqref{eq:rel5}
generate the ideal $I_R$.
Hence \eqref{eq:FellMnscalarrel},\eqref{eq:rel1-4},\eqref{eq:rel5}
consitute the defining relations of the algebra $\Fell(M(n))$.
\end{prop}
\begin{proof}
Assume we have a relation of the form \eqref{eq:RLLmatrixelements_residual_general}
and that a limit in one of the terms, $\lim_{w\to z}\ph(w)R_{xy}^{ab}(\la,w)$, say,
exists and is nonzero.
Then one of the following cases occurs.
\begin{enumerate}
\item At $w=z$, $\ph(w)$ and $R_{xy}^{ab}(\la,w)$ are both regular. If this holds for all terms,
then the relation is just a multiple of one of \eqref{eq:rel1}-\eqref{eq:rel4}.
\item At $w=z$, $\ph(w)$ has a pole while $R_{xy}^{ab}(\la,w)$ is regular.
Then $R_{xy}^{ab}(\la,w)$ must vanish identically at $w=z$.
The only case where this is possible is when $x\neq y$ and
$R_{xy}^{ab}(\la,w)=\al(\la_{xy},w)$ and $z=p^k$. But then there is another
term containing $\be$ which is never identically zero for any $z$,
and hence the limit in that term does not exist.
\item At $w=z$, $\ph(w)$ is regular while $R_{xy}^{ab}(\la,w)$ has a pole.
Since these poles are simple and occur only when $z\in q^{-2}p^\Z$,
the function $\ph$ must have a zero of multiplicity one there.
We can assume without loss of generality that
$\ph$ has the specific form
 \[\ph(w)=\frac{q^{-1}\th(q^2w)}{q\th(q^{-2}w)}.\]
Then, if $a\neq c$ and $b\neq d$, 
\eqref{eq:RLLmatrixelements_residual_general} becomes 
the residual RLL relation \eqref{eq:rel5}.
 
If instead $c=a$, $b\neq d$, and we take
$z_1=z$, $z_2=p^kq^2z$ in \eqref{eq:RLLmatrixelements_residual_general}
we get, using \eqref{eq:albelimits},
\[
0=\al(\rh_{bd},p^kq^2)e_{ad}(p^kq^2z)e_{ab}(z)-
\be(\rh_{bd},p^kq^2)e_{ab}(p^kq^2z)e_{ad}(z),
\]
or, rewritten,
\[
e_{ad}(p^kq^2z)e_{ab}(z)=q^{2k\rh_{bd}}\frac{E(\rh_{bd}-1)}{E(\rh_{bd}+1)}
e_{ab}(p^kq^2z)e_{ad}(z).
\]
However this relation is already derivable from \eqref{eq:rel2} as follows.
Take $z_1=p^kq^2z$ and $z_2=z$
in \eqref{eq:rel2} and multiply both sides by
$q^{2k\rh_{bd}}\frac{E(\rh_{bd}-1)}{E(\rh_{bd}+1)}$ and then
use \eqref{eq:rel2} on the right hand side.

Similarly, if $a\neq c$, $d=b$,
$z_1=z$, $z_2=p^kq^2z$, $\ph(w)=\frac{q^{-1}\th(q^2w)}{q\th(q^{-2}w)}$
in \eqref{eq:RLLmatrixelements_residual_general}
and using \eqref{eq:albelimits} we get
\[
\al(\la_{ac},p^kq^2)e_{ab}(z)e_{cb}(p^kq^2z)-\be(\la_{ca},p^kq^2)e_{cb}(z)e_{ab}(p^kq^2z)=0,
\]
or,
\[
e_{ab}(z)e_{cb}(p^kq^2z)=q^{2k\la_{ca}}e_{cb}(z)e_{ab}(p^kq^2z).
\]
Similarly to the previous case, this identity follows already from \eqref{eq:rel3}.

\end{enumerate}

\end{proof}

\section{Left and right elliptic exterior algebras} \label{sec:extalgs}
\subsection{Corepresentations of $\hf$-bialgebroids}
We recall the definition of corepresentations of an $\hf$-bialgebroid
given in \cite{KoRo}.
\begin{dfn}
An \emph{$\hf$-space} $V$ is an $\hf^*$-graded vector space over $M_{\hf^*}$,
$V=\bigoplus_{\al\in\hf^*} V_\al$,
where each $V_\al$ is $M_{\hf^*}$-invariant. A morphism of $\hf$-spaces is a degree-preserving
$M_{\hf^*}$-linear map.
\end{dfn}
Given an $\hf$-space $V$ and an $\hf$-bialgebroid $A$, we define
$A\widetilde{\otimes} V$ to be the $\hf^*$-graded space with
$(A\widetilde{\otimes} V)_\al = \bigoplus_{\be\in\hf^*} (A_{\al\be}\otimes_{M_{\hf^*}} V_\be)$
where $\otimes_{M_{\hf^*}}$ denotes $\otimes_\C$ modulo the relations
\[\mu_r(f)a \otimes v = a\otimes f v,\]
for $f\in M_{\hf^*}$, $a\in A$, $v\in V$.
$A\widetilde{\otimes}V$ becomes an $\hf$-space with the $M_{\hf^*}$-action
$f (a\otimes v)= \mu_l(f)a\otimes v$.
Similarly we define $V\widetilde{\otimes} A$ as an $\hf$-space by
$(V\widetilde{\otimes} A)_\be = \bigoplus_\al V_\al\otimes_{M_{\hf^*}} A_{\al\be}$
where $\otimes_{M_{\hf^*}}$ here means $\otimes_\C$ modulo the relation
$v \otimes \mu_l(f)a =fv\otimes a$, and $M_{\hf^*}$-action given by
$f(v\otimes a)=v\otimes \mu_r(f)a$.

For any $\hf$-space $V$ we have isomorphisms $D_\hf\widetilde{\otimes} V\simeq V\simeq V\widetilde{\otimes} D_\hf$ given by
\begin{equation}\label{eq:tensorunit2}
T_{-\al}\otimes v \simeq v \simeq v\otimes T_{\al},\quad \text{for $v\in V_\al$},
\end{equation}
extended to $\hf$-space morphisms.

\begin{dfn}
A left corepresentation $V$ of an $\hf$-bialgebroid $A$ is an $\hf$-space equipped with
an $\hf$-space morphism $\De_V:V\to A\widetilde{\otimes} V$ 
such that $(\De_V\otimes 1)\circ\De_V = (1\otimes \De)\circ\De_V$
and $(\ep\otimes 1)\circ \De_V = \Id_V$ (under the identification \eqref{eq:tensorunit2}).
\end{dfn}
\begin{dfn} A left $\hf$-comodule algebra $V$ over an $\hf$-bialgebroid $A$
is a left corepresentation
$\De_V:V\to A\widetilde{\otimes} V$ and in addition a $\C$-algebra such that
$V_\al V_\be\subseteq V_{\al+\be}$ and such that $\De_V$ is an algebra morphism,
when $A\widetilde{\otimes} V$ is given the natural algebra structure.
\end{dfn}
Right corepresentations and comodule algebras are defined analogously.

\subsection{The comodule algebras $\La$ and $\La'$.}
We define in this section an elliptic analog of the exterior algebra,
following \cite{KovN}, where it was carried out in the
trigonometric non-spectral case. It will lead to natural definitions
of elliptic minors as certain elements of $\Fell(M(n))$.
One difference between this approach and the one in \cite{FeVa97}
is that the elliptic exterior algebra in our setting is
really an algebra, and not just a vector space. Another one is that
the commutation relations in our elliptic exterior algebras 
are completely determined by requiring the natural relations
\eqref{eq:leftextalgrel1},
\eqref{eq:leftextalgrel2}, \eqref{eq:extalgcoproddef},
and that the coaction is an algebra homomorphism.
This fact can be seen from the proof of Proposition \ref{prop:extalgfirst}.
Since the proof does not depend on the particular form of $\al$ and $\be$,
we can obtain exterior algebras for any $\hf$-bialgebroid
obtain through the generalized FRST-construction from
an R-matrix in the same manner.
In particular the method is independent of the gauge equivalence class
of $R$.

Let $\La$ be the unital associative $\C$-algebra generated by $v_i(z)$,
$1\le i\le n$, $z\in \C^\times$ and a copy of $M_{\hf^*}$ embedded as a
subalgebra subject to the relations
\begin{subequations}
\begin{gather}
\label{eq:leftextalgrel1}
f(\ze)v_i(z)=v_i(z)f(\ze+\om(i)),\\
\label{eq:leftextalgrel2}
v_i(z)v_i(w)=0,\\
\label{eq:leftextalgrel3}
\al(\ze_{kj},z/w)v_k(z)v_j(w)+\be(\ze_{kj},z/w)v_j(z)v_k(w)=0,
\end{gather}
for $i,j,k\in\oneton$, $j\neq k$, $z,w\in\C^\times$, $z/w\notin\{p^kq^{-2}:k\in\Z\}$
and $f\in M_{\hf^*}$.
We require also the residual relation of
\eqref{eq:leftextalgrel3} obtained by multiplying by
$\ph(z/w)=\frac{q^{-1}\th(q^2z/w)}{q\th(q^{-2}z/w)}$ and letting $z/w\to p^{-k}q^{-2}$.
After simplification using \eqref{eq:albelimits}, we get
\begin{equation}\label{eq:leftextalgrel4}
 v_k(z)v_j(p^kq^2z)=q^{2k\ze_{jk}}v_j(z)v_k(p^kq^2z).
\end{equation}
\end{subequations}

$\La$ becomes an $\hf$-space by
\[\mu_\La(f)v=f(\ze)v\]
and requiring $v_i(z)\in \La_{\om(i)}$ for each $i,z$.

\begin{prop} \label{prop:extalgfirst}
 $\La$ is a left comodule algebra over $\Fell(M(n))$ with
left coaction $\De_\La :\La\to \Fell(M(n))\widetilde{\otimes}\La$ satisfying
\begin{align}
\label{eq:extalgcoproddef}
\De_\La(v_i(z)) &= \sum_j e_{ij}(z)\otimes v_j(z),\\
\De_\La(f(\ze)) &= f(\la)\otimes 1.
\end{align}
\end{prop}
\begin{proof}
We have
\begin{align*}
&\De_\La(v_i(z))\De_\La(v_i(w))=\sum_{jk} e_{ij}(z)e_{ik}(w)\otimes v_j(z)v_k(w)=\\
&=\sum_{j\neq k}\big(\al(\mu_{jk},\frac{z}{w})e_{ik}(w)e_{ij}(z)+\be(\mu_{kj},\frac{z}{w})e_{ij}(w)e_{ik}(z)
\big)\otimes v_j(z)v_k(w)=\\
&=\sum_{j\neq k} e_{ij}(w)e_{ik}(z)\otimes\big(\al(\ze_{kj},\frac{z}{w})v_k(z)v_j(w)+
\be(\ze_{kj},\frac{z}{w})v_j(z)v_k(w)\big)=0.
\end{align*}
Similarly one proves that \eqref{eq:leftextalgrel3},\eqref{eq:leftextalgrel4} are preserved.
\end{proof}

Relation \eqref{eq:leftextalgrel3} is not symmetric under
interchange of $j$ and $k$. We now
derive a more explicit, independent, set of relations for $\La$.
We will use the function $E$, defined in \eqref{eq:Edef}.

\begin{prop}\label{prop:extalg}
(i) The following is a complete set of relations for $\La$
\begin{subequations}\label{eq:leftsimp0-2}
\begin{align}
 \label{eq:leftsimp0}
& f(\ze)v_i(z)=v_i(z)f(\ze+\om(i)),\\
 \label{eq:leftsimp1}
&v_k(p^sq^2z)v_j(z)=-q^{2s\ze_{kj}}\frac{E(\ze_{kj}-1)}{E(\ze_{kj}+1)}v_j(p^sq^2z)v_k(z),\quad\forall s\in\Z, k\neq j,\\
\label{eq:leftsimp1.5}
&v_k(z)v_j(p^sq^2z)=q^{2s\ze_{jk}}v_j(z)v_k(p^sq^2z),\\
\label{eq:leftsimp2}
&v_k(z)v_j(w)=0\quad\text{if $z/w\notin\{p^sq^{\pm 2}|s\in\Z\}$ or if $k=j$.}
\end{align}
\end{subequations}

(ii)
The set
\begin{equation}\label{eq:leftextalgbasis}
 \{v_{i_d}(z_d)\cdots v_{i_1}(z_1) : 1\le i_1<\cdots <i_d\le n,
 \frac{z_{i+1}}{z_i}\in p^\Z q^{\pm 2}\}
 \end{equation}

is a basis for $\La$ over $M_{\hf^*}$.
 
\end{prop}

\begin{proof}
(i) Elimination of the $v_j(z)v_k(w)$-term in \eqref{eq:leftextalgrel3} yields
\begin{equation}\label{eq:leftel1}
\big(\al(\ze_{jk},\frac{z}{w})\al(\ze_{kj},\frac{z}{w})-\be(\ze_{kj},\frac{z}{w})\be(\ze_{jk},\frac{z}{w})\big)
v_k(z)v_j(w)=0.
\end{equation}
Combining \eqref{eq:leftel1}, \eqref{eq:albeth} and the fact that the $\th(z)$ is zero
precisely for $z\in\{p^k|k\in\Z\}$ we deduce that in $\La$,
\begin{equation}\label{eq:leftel2}
v_k(z)v_j(w)\neq 0 \Longrightarrow \frac{z}{w}=p^sq^2 \text{ for some } s\in\Z.
\end{equation}
Using \eqref{eq:albepk} we obtain from \eqref{eq:leftel2} and \eqref{eq:leftextalgrel2},\eqref{eq:leftextalgrel3}
that relations \eqref{eq:leftsimp1},\eqref{eq:leftsimp2} hold
in the left elliptic exterior algebra $\La$.
Relations \eqref{eq:leftsimp0},\eqref{eq:leftsimp1.5}
are just repetitions of \eqref{eq:leftextalgrel1},\eqref{eq:leftextalgrel4}.

(ii) It follows from the relations that each monomial in $\La$ can be uniquely written
as $f(\ze)v_{i_d}(z_d)\cdots v_{i_1}(z_1)$ where $1\le i_1<\cdots<i_d\le n$ and $f\in M_{\hf^*}$.
It remains to show that the set \eqref{eq:leftextalgbasis} is linearly independent over $M_{\hf^*}$.
Assume that a linear combination of basis elements is zero, and that the sum has minimal number of terms.
By multiplying from the right or left by $v_j(w)$ for appropriate $j$, $w$
we can assume the sum is of the form
\[f_1(\ze)v_{i_d}(z_d^1)\cdots v_{i_1}(z_1^1)+\cdots+
f_r(\ze)v_{i_d}(z_d^r)\cdots v_{i_1}(z_1^r) =0 \]
for some fixed set $\{i_1,\ldots,i_d\}$. By the relations, a monomial $v_{i_d}(z_d)\cdots v_{i_1}(z_1)$
can be given the ``degree'' $\sum_{i=1}^d z_it^{i-1}\in \C[t]$, where $t$ is an indeterminate.
Formally, consider
$\C(t)\otimes \La$,
the tensor product (over $\C$) of $\La$ by the field of rational functions in $t$.
We identify $\La$ with its image under $\La\ni v\mapsto 1\otimes v\in \C(t)\otimes \La$,
and view $\C(t)\otimes\La$ naturally as a vector space over $\C(t)$.
By relations \eqref{eq:leftsimp0}-\eqref{eq:leftsimp2}, there is a $\C$-algebra automorphism
$T$ of $\C(t)\otimes \La$ satisfying
$T(v_j(z))=t v_j(z)$, $T(f(\ze))=f(\ze)$ and $T(p\otimes 1)=p\otimes 1$.
Define
\[D(v_i(z))=z v_i(z),\quad D(f(\ze))=0,\quad D(p\otimes 1)=0,\]
for $f\in M_{\hf^*}$, $p\in \C(t)$ and $i\in\oneton$, $z\in\C^\times$
and extend $D$ to a $\C$-linear map $D:\C(t)\otimes\La\to \C(t)\otimes \La$ by
requiring
\begin{equation}\label{eq:skewderivation}
D(ab)=D(a)T(b)+aD(b)
\end{equation}
for $a,b\in \C(t)\otimes \La$.
The point is that the requirement \eqref{eq:skewderivation}
respects relations \eqref{eq:leftsimp0}-\eqref{eq:leftsimp2}, making $D$ well defined.
Write $u_j=f_j(\ze)v_{i_d}(z_d^j)\cdots v_{i_1}(z_1^j)$.
Then one checks that $D(u_j)=p_j(t) u_j$, where $p_j(t)=\sum_{i=1}^d z_i^j t^{i-1}$.
By applying $D$ repeatedly we get
\begin{align*}
u_1(z^1)+&\cdots +u_r(z^r)=0,\\
p_1(t)u_1(z^1)+&\cdots +p_r(t)u_r(z^r)=0,\\
& \;\; \vdots\\
p_1(t)^{r-1}u_1(z^1)+&\cdots +p_r(t)^{r-1}u_r(z^r)=0.
\end{align*}
Inverting the Vandermonde matrix $(p_j(t)^{i-1})_{ij}$ we obtain $u_j(z^j)=0$ for each $j$,
i.e. $f_j(\ze)=0$ for each $j$. This proves linear independence of \eqref{eq:leftextalgbasis}.
\end{proof}

Analogously one defines a right comodule algebra $\La'$ with generators
$w^i(z)$ and $f(\ze)\in M_{\hf^*}$. The following relations will be used:
\begin{subequations}
\begin{align}
\label{eq:rightsimp1}
&w^k(z)w^j(p^sq^2z)=-q^{2s\ze_{kj}}w^j(z)w^k(p^sq^2z),\quad\forall s\in\Z, k\neq j,\\
\label{eq:rightsimp2}
&w^k(z_1)w^j(z_2)=0\quad\text{if $z_2/z_1\notin\{p^sq^{\pm 2}|s\in\Z\}$ or if $k=j$.}
\end{align}
\end{subequations}
$\La'$ has also $M_{\hf^*}$-basis of the form \eqref{eq:leftextalgbasis}. In fact
$\La$ and $\La'$ are isomorphic as algebras.

\subsection{Action of the symmetric group}\label{sec:Snaction}
From \eqref{eq:rel1-4},\eqref{eq:rel5} we see that
$S_n\times S_n$ acts by $\C$-algebra automorphisms on $\Fell(M(n))$ as follows
\[(\si,\tau)(f(\la))=f(\la\circ L_\si), \qquad (\si,\tau)(f(\mu))=f(\mu\circ L_\ta),\]
\[(\si,\tau)(e_{ij}(z))=e_{\si(i)\tau(j)}(z),\]
where $L_\si :\hf\to \hf$ ($\si\in S_n$) is given by permutation of coordinates:
\[L_\si(\diag(h_1,\ldots,h_n))=\diag(h_{\si(1)},\ldots,h_{\si(n)}).\]
Also, $S_n$ acts on $\La$ by $\C$-algebra automorphisms via
\begin{equation}\label{eq:SnonLa}
\si(f(\ze)) = f(\ze\circ L_\si),\qquad
\si(v_i(z))=v_{\si(i)}(z).
\end{equation}
Similarly we define an $S_n$ action on $\La'$.
\begin{lem} For each $v\in\La$, $w\in\La'$ and any $\si,\ta\in S_n$ we have
\begin{align} \label{eq:leftSncovariance}
\De_\La(\si(v)) &= ((\si,\ta)\otimes \ta )(\De_\La(v)),\\
\label{eq:rightSncovariance}
\De_{\La'}(\ta(w)) &= (\si\otimes (\si,\ta))(\De_{\La'}(w)).
\end{align}
\end{lem}
\begin{proof}
By multiplicativity, it is enough to prove these claims on the generators,
which is easy.
\end{proof}

\section{Elliptic quantum minors} \label{sec:minors}
\subsection{Definition}
For $I\subseteq\oneton$ we set
\begin{align}
\label{eq:FIdef}F_I(\ze)&=\prod_{i,j\in I, i<j} E(\ze_{ij}+1),&
F^I(\ze)&=\prod_{i,j\in I, i<j} E(\ze_{ij}),
\end{align}
and define the left and right elliptic sign functions
\begin{align}\label{eq:leftsigndef}
\sgn_I(\si;\ze) &=\frac{\si(F_I(\ze))}{F_{\si(I)}(\ze)}= \prod_{i,j\in I, i<j, \si(i)>\si(j)}
 \frac{E(\ze_{\si(i)\si(j)}+1)}{E(\ze_{\si(j)\si(i)}+1)},\\
\label{eq:rightsigndef}
\sgn^I(\si;\ze) &=\frac{F^{\si(I)}(\ze)}{\si(F^I(\ze))}= \prod_{i,j\in I, i<j, \si(i)>\si(j)}
 \frac{E(\ze_{\si(j)\si(i)})}{E(\ze_{\si(i)\si(j)})},
\end{align}
for $\si\in S_n$.
In fact, $E(\ze_{ij})/E(\ze_{ji})=-1$ so $\sgn^{\oneton}(\si;\ze)$ is just
the usual sign $\sgn(\si)$. However we view this as a ``coincidence''
depending 
on the particular choice of R-matrix from its gauge equivalence class.
We keep our notation to emphasize that the methods do not depend on
this choice of $R$-matrix.

We will denote the elements of a subset $I\subseteq\oneton$
by $i_1<i_2<\cdots$.
\begin{prop}
Let $I\subseteq\oneton$, d=\#I, $\si\in S_n$ and $J=\si(I)$.
Then for $z\in\C^\times$,
\begin{equation} \label{eq:minors1}
v_{\si(i_d)}(q^{2(d-1)}z)\cdots v_{\si(i_1)}(z)=
\sgn_I(\si;\ze)
v_{j_d}(q^{2(d-1)}z)\cdots v_{j_1}(z)
\end{equation}
and
\begin{equation} \label{eq:minors2}
w^{\si(i_1)}(z)\cdots w^{\si(i_d)}(q^{2(d-1)}z)=
\sgn^I(\si;\ze)w^{j_1}(z)\cdots w^{j_d}(q^{2(d-1)}z).
\end{equation}
\end{prop}
\begin{proof}
We prove \eqref{eq:minors1}. The proof of \eqref{eq:minors2} is analogous.
We proceed by induction on $\#I=d$, the case $d=1$ being clear.
If $d>1$, set $I'=\{i_1,\ldots,i_{d-1}\}, J'=\si(I')$.
Let $1\le j_1'<\cdots <j_{d-1}'\le n$ be the elements of $J'$.
By the induction hypothesis, the left hand side of \eqref{eq:minors1}
equals
\begin{equation}
v_{\si(i_d)}(q^{2(d-1)}z)\sgn_{I'}(\si,\ze) v_{j_{d-1}'}(q^{2(d-2)}z)\cdots v_{j_1'}(z).
\end{equation}
Now $v_{\si(i_d)}(q^{2(d-1)}z)$ commutes with $\sgn_{I'}(\si,\ze)$
since the latter only involves $\ze_{ij}$ with $i,j\neq \si(i_d)$. Using
the commutation relations \eqref{eq:leftsimp1} we obtain
\begin{equation}
\sgn_{I'}(\si,\ze)\cdot\prod_{j\in J', j>\si(i_d) }\frac{E(\ze_{j\si(i_d)}+1)}{E(\ze_{\si(i_d)j}+1)} \cdot
v_{j_d}(q^{2(d-1)}z)\cdots v_{j_1}(z).
\end{equation}
Replace $j\in J'$ such that $j>\si(i_d)$ by $\si(i)$ where $i\in I, i<i_d, \si(i)>\si(i_d)$.
\end{proof}

Introduce the normalized monomials
\begin{align}
 \label{eq:vIdef}
 v_I(z)&=F_I(\ze)^{-1} v_{i_r}(q^{2(d-1)}z)v_{i_{r-1}}(q^{2(d-2)}z)\cdots v_{i_1}(z) \in\Lambda,\\
 \label{eq:wIdef}
 w^I(z) &= F^I(\ze) w^{i_1}(z)w^{i_2}(q^2z)\cdots w^{i_d}(q^{2(d-1)}z) \in \Lambda'.
\end{align}

\begin{cor} \label{cor:Snonvw}
Let $I\subseteq\oneton$. For any permutation $\si\in S_n$,
\begin{align}
\si(v_I(z))&=v_{\si(I)}(z),&\si(w^I(z))&=w^{\si(I)}(z),
\end{align}
for any $z\in\C^\times$. In particular $v_I(z)$ and $w^I(z)$
are fixed by any permutation which preserves the subset $I$.
\end{cor}
\begin{proof} Let $J=\si(I)$. Then
\begin{align*}
\si(v_I(z))&=\si(F_I(\ze)^{-1})v_{\si(i_d)}(q^{2(d-1)}z)\cdots v_{\si(i_1)}(z)=\\
&=\si(F_I(\ze))^{-1}\sgn_I(\si;\ze) v_{j_d}(q^{2(d-1)}z)\cdots v_{j_1}(z)=v_{\si(I)}(z).
\end{align*}
The proof for $w^I(z)$ is analogous.
\end{proof}

We are now ready to define certain elements of the $\hf$-bialgebroid
$\Fell(M(n))$ which are analogs of minors.
\begin{prop}
For $I,J\subseteq\oneton$ and $z\in\C^\times$,
the left and right elliptic minors, $\overleftarrow{\xi}_I^J(z)$ and
$\overrightarrow{\xi}_I^J(z)$ respectively, can be defined by
\begin{align} \label{eq:minors3}
 \De_\La(v_I(z))&=\sum_J
 \overleftarrow{\xi}_I^J(z)\otimes v_J(z),\\
 \label{eq:minors4}
 \De_{\La'}(w^J(z))&=\sum_I w^I(z)\otimes \overrightarrow{\xi}_I^J(z),
\end{align}
where the sums are taken over all subsets of $\oneton$.

If $\#I\neq\#J$, then $\overleftarrow{\xi}_I^J(z)=0=\overrightarrow{\xi}_I^J(z)$ for all
$z$. If $\#I=\#J=d$, they are explicitly given by
\begin{multline}\label{eq:leftminor}
\overleftarrow{\xi}_I^J(z)= \\
 =\frac{F_J(\rh)}{F_I(\la)} \sum_{\ta\in S_J} \frac{\sgn_J(\ta;\rh)}{\sgn_I(\si;\la)}
 e_{\si(i_d)\ta(j_d)}(q^{2(d-1)}z)
 e_{\si(i_{d-1})\ta(j_{d-1})}(q^{2(d-2)}z)\cdots e_{\si(i_1)\ta(j_1)}(z)
\end{multline}
for any $\si\in S_I$, and
\begin{equation}\label{eq:rightminor}
\overrightarrow{\xi}_I^J(z)= \frac{F^J(\rh)}{F^I(\la)} \sum_{\si\in S_I} \frac{\sgn^J(\ta;\rh)}{\sgn^I(\si;\la)}
 e_{\si(i_1)\ta(j_1)}(z)e_{\si(i_2)\ta(j_2)}(q^2z)\cdots e_{\si(i_d)\ta(j_d)}(q^{2(d-1)}z)
\end{equation}
for any $\ta\in S_J$.
Moreover, 
\begin{equation}\label{eq:Snonminors}
(\si,\ta)(\overleftarrow{\xi}_I^J(z))=\overleftarrow{\xi}_{\si(I)}^{\ta(J)}(z)
 \quad\text{and}\quad (\si,\ta)(\overrightarrow{\xi}_I^J(z))=
 \overrightarrow{\xi}_{\si(I)}^{\ta(J)}(z)
\end{equation}
for any $(\si,\tau)\in S_n\times S_n$ and $z\in\C^\times$.
\end{prop}
\begin{rem}
In Theorem \ref{thm:leftequalsright}
we will prove that, in fact, $\overleftarrow{\xi}_I^J(z)=\overrightarrow{\xi}_I^J(z)$.
\end{rem}
\begin{proof}
We prove the statements concerning the left elliptic minor $\overleftarrow{\xi}_I^J(z)$.
We have
\begin{align*}
&\De_\La(v_I(z))=\sum_{1\le k_1,\ldots,k_d\le n} F_I(\la)^{-1}e_{i_1k_1}(q^{2(d-1)}z)\cdots
 e_{i_dk_d}(z)\otimes v_{k_1}(q^{2(d-1)}z)\cdots v_{k_d}(z)=\\
&=\sum_{J, \#J=d} \sum_{\ta\in S_J} F_I(\la)^{-1}
e_{i_1\ta(j_1)}(q^{2(d-1)}z)\cdots e_{i_d\ta(j_d)}(z)\otimes
v_{\ta(j_1)}(q^{2(d-1)}z)\cdots v_{\ta(j_d)}(z)=\\
&=\sum_{J, \#J=d} \left( \sum_{\ta\in S_J} \frac{\ta(F_J(\rh))}{F_I(\la)}
e_{i_1\ta(j_1)}(q^{2(d-1)}z)\cdots e_{i_d\ta(j_d)}(z)\right) \otimes v_J(z).
\end{align*}
Thus \eqref{eq:minors3} holds when $\overleftarrow{\xi}_I^J(z)$ is defined by
\eqref{eq:leftminor} with $\si=\Id$. Then the right hand side of \eqref{eq:leftminor}
equals $(\si,\Id)(\overleftarrow{\xi}_I^J(z))$. Thus only \eqref{eq:Snonminors} remains.
Using \eqref{eq:leftSncovariance} and Corollary \ref{cor:Snonvw} we have
\[\De_\La\big(\si(v_I(z))\big)=((\si,\ta)\otimes \ta)\left(\De_\La(v_I(z))\right)
=\sum_J (\si,\ta)(\overleftarrow{\xi}_I^J(z))\otimes v_{\ta(J)}(z)).\]
On the other hand, again by Corollary \ref{cor:Snonvw},
\[\De_\La\big(\si(v_I(z))\big)=\De_\La(v_{\si(I)}(z))
=\sum_J \overleftarrow{\xi}_{\si(I)}^{\ta(J)}(z)\otimes v_{\ta(J)}(z)),\]
where we made the substitution $J\mapsto \ta(J)$.
This proves the first equality in \eqref{eq:Snonminors}.
%
The statements concerning the right elliptic minors are proved analogously.
\end{proof}

\subsection{Equality of left and right minors} \label{sec:cherednik}
The goal of this section is to prove Theorem \ref{thm:leftequalsright}
stating that the left and right elliptic minors coincide.
We use ideas from Section 3 of \cite{FeVa97}, where the authors
study the objects of Felder's tensor category \cite{Fe95a}
and associate a linear operator (product of R-matrices) on $V^{\otimes n}$
to each diagram of a certain form,
a kind of braid group representation. Then they consider the operator
associated to the longest permutation, in \cite{ZhShYu} called the Cherednik operator.
Instead of working with representations,
we proceed inside the $\hf$-bialgebroid $\Fell(M(n))$ and
consider certain operators on $V^{\otimes n}\otimes \Fell(M(n))$
depending on $n$ spectral parameters.
Using the analog of the Cherednik operator
we prove an extended RLL-relation \eqref{eq:WLLrel}.
Theorem \ref{thm:leftequalsright} then follows by extracting matrix
elements from both sides of this matrix equation.

In this section, we set $\Fscr=\Fell(M(n))$. 
Recall the operators
from Section \ref{sec:operatorform}, defined for any $\hf$-bialgebroid
$A_R$ obtained from the FRST-construction.
When specializing to $\Fscr$ we get operators $\R(\la,z)$, $\R(\rh,z)\in\End(V\otimes V\otimes\Fscr)$,
where $V=\C^n$.
For $z\in\C^\times$, define the following linear operators on $V^{\otimes n}\otimes \Fscr$:
\[\R^{ij}(\la,z) := \lim_{w\to z}\th(q^2w)\R(\la,w)^{(i,j,n+1)},
 \qquad \R^{ij}(\rh,z) :=\lim_{w\to z} \th(q^2w)\R(\rh,w)^{(i,j,n+1)}.\]
The upper indices in parenthesis are tensor leg numbering and indicate the tensor factors
the operator should act on.
The limits are taken in the sense of taking limits of each matrix element.
These operators are well-defined for any $z$, since we
multiply away the singularities in $z$ of $\al$ and $\be$ \eqref{eq:aldef},\eqref{eq:bedef}.

Let $\mathcal{E}_n$ denote the algebra of all functions
\[F:(\C^\times)^n\to \End(V^{\otimes n}\otimes \Fscr).\]
%
The symmetric group $S_n$ acts on $\mathcal{E}_n$ by
\begin{equation}\label{eq:38924}
\si(F(z))=(\si\otimes \Id_\Fscr)\circ F(\si(z)) \circ (\si^{-1}\otimes \Id_\Fscr)
\end{equation}
for $F(z)\in\mathcal{E}_n$ and $\si\in S_n$.
In the right hand side of \eqref{eq:38924},
$\si$ acts on $(\C^\times)^n$ by permuting coordinates, and on $V^{\otimes n}$
by permuting the tensor factors.
For example we have
\[(23)\left(\R^{12}(\la,z_1/z_2)\right)=\R^{13}(\la,z_1/z_3).\]
Consider the skew group algebra $\mathcal{E}_n \ast S_n$,
defined as the algebra with underlying space 
$\mathcal{E}_n\otimes \C S_n$, 
where $\C S_n$ is the group algebra, with the multiplication
\begin{equation}\label{eq:skewproduct}
 (F(z)\otimes\si)(G(z)\otimes\ta) = F(z)\si(G(z))\otimes \si\ta
\end{equation}
for $\si,\ta\in S_n$, $F(z),G(z)\in \mathcal{E}_n$. Since $\si$ acts on $\mathcal{E}_n$ by automorphisms,
$\mathcal{E}_n\ast S_n$ is an associative algebra. The constant function
$z\mapsto \Id_{V^{\otimes n}\otimes \Fscr}\otimes (1)$ is the unit element.
Let $B_n$ be the monoid (set with unital associative multiplication)
generated by $\{s_1,\ldots,s_{n-1}\}$ and relations
\begin{align*}
s_is_{i+1}s_i&=s_{i+1}s_is_{i+1}\quad\text{for $1\le i\le n-2$,}\\
 s_is_j&=s_js_i\quad\text{if $|i-j|>1$}.
\end{align*}
Let $\si_i=(i\;\; i+1)\in S_n$.
We have an epimorphism $\pi:B_n\to S_n$ given by $\pi(s_i)=\si_i$, $\pi(1)=(1)$.
Define
\[\Cdd(1) = \Id_{V^{\otimes n}\otimes \Fscr}\otimes (1),\]
\[\Cdd(s_i) = \R^{i,i+1}(\la-h^{\ge i+2},z_i/z_{i+1})\otimes \si_i.\]
Here and below we use $h^{\ge k}$ to denote the expression $\sum_{j=k}^n h^j$.

\begin{prop}
$\Cdd$ extends to a well-defined morphism of monoids, i.e. a map
\[\Cdd: B_n \to \mathcal{E}_n \ast S_n\]
satisfying $\Cdd(b_1b_2)=\Cdd(b_1)\Cdd(b_2)$ for any $b_1,b_2\in B_n$.
\end{prop}
\begin{proof}
We have to show the relations
\begin{align}
 \label{eq:doprel2}\Cdd(s_i)\Cdd(s_{i+1})\Cdd(s_i)&=
 \Cdd(s_{i+1})\Cdd(s_i)\Cdd(s_{i+1}),\\
 \label{eq:doprel3}\Cdd(s_i)\Cdd(s_j)&=\Cdd(s_j)\Cdd(s_i)\quad\text{if $|i-j|>1$}.
\end{align}
Relation \eqref{eq:doprel2} follows from the QDYBE \eqref{eq:qdybe}.
For example, $\Cdd(s_i)\Cdd(s_{i+1})\Cdd(s_i)$ equals
\[\R^{i,i+1}(\la-h^{\ge i+2},\frac{z_i}{z_{i+1}})
\R^{i,i+2}(\la-h^{\ge i+3},\frac{z_i}{z_{i+2}})
\R^{i+1,i+2}(\la-h^i-h^{\ge i+3},\frac{z_{i+1}}{z_{i+2}})
\otimes \si_i\si_{i+1}\si_i\]
Relation \eqref{eq:doprel3} is easy to check, using the $\hf$-invariance of $R$.
\end{proof}

For $b\in B_n$ we define
$\Che_b(\la,z)\in \mathcal{E}_n$ by 
\begin{equation}\label{eq:Chedef}
\Cdd(b) = \Che_b(\la,z)\otimes \pi(b).
\end{equation}
From this and the product rule \eqref{eq:skewproduct} follows that
\begin{equation}\label{eq:chemultiplicativity}
\Che_{b_1b_2}(\la,z)=\Che_{b_1}(\la,z)\cdot \pi(b_1)\left(\Che_{b_2}(\la,z)\right)
\end{equation}
for $b_1,b_2\in B_n$.
By replacing $\la$ by $\rh$ we get similarly operators $\Che_b(\rh,z)$.

Recall the operators $\L(z)\in\End(V\otimes\Fscr)$
 from Section \ref{sec:operatorform}.
Define for $z\in\C^\times$, $i\in\oneton$,
\[\L^i(z)=\L(z)^{(i,n+1)}\in \End(V^{\otimes n}\otimes \Fscr).\]

If $i,j,k$ are distinct, then one can check that
\begin{align}
\label{eq:RLLR1} \R^{ij}(\la-h^k,z)\L^k(w)&=\L^k(w)\R^{ij}(\la,z),\\
\label{eq:RLLR2} \R^{ij}(\rh,z)\L^k(w)&=\L^k(w)\R^{ij}(\rh+h^k,z).
\end{align}

Due to the RLL relations \eqref{eq:RLLmatrixelements_residual_general} 
we have
\begin{equation} \label{eq:RLLwithmodifiedR}
 \R^{12}(\la,\frac{z_1}{z_2})\L^1(z_1)\L^2(z_2)=\L^2(z_2)\L^1(z_1)\R^{12}(\rh+h^1+h^2,\frac{z_1}{z_2})
\end{equation}
for any $z_1,z_2\in\C^\times$.

Define $t_d\in B_n$, $d\in\oneton$, recursively by
\[t_d= \begin{cases}
 t_{d-1} s_{d-1}s_{d-2}\cdots s_1, & d>1\\
 1,&d=1.
\end{cases}
\]
Let $\ta_d$ be the image of $t_d$ in $S_n$:
\[\ta_d:=\pi(t_d)=
\mbox{$
\begin{pmatrix}1&2&\cdots &d &d+1&\cdots &n\\
d&d-1&\cdots & 1 & d+1 &\cdots & n\end{pmatrix}\in S_n$.}\]

\begin{prop} Let $1\le d\le n$. For any $z=(z_1,\ldots,z_d)\in (\C^\times)^d$ we have
 \begin{equation} \label{eq:WLLrel}
 \Che_{t_d}(\la,z)\L^1(z_1)\cdots\L^d(z_d)=
 \L^d(z_d)\cdots \L^1(z_1)\Che_{t_d}(\rh+h^{\le d},z).
\end{equation}
\end{prop}
\begin{proof}
We use induction on $d$. The case $d=1$ is trivial, while
$d=2$ is the RLL relation \eqref{eq:RLLwithmodifiedR}.
If $d>2$, write $t_d =t_{d-1}u_d$, where $u_d=s_{d-1}s_{d-2}\cdots s_1$.
Thus, by \eqref{eq:chemultiplicativity},
\begin{equation}\label{eq:chefactorization}
\Che_{t_d}(\la,z)=\Che_{t_{d-1}}(\la,z)\cdot \ta_{d-1}\left(\Che_{u_d}(\la,z)\right).
\end{equation}
We claim that
\begin{multline} \label{eq:CLLclaim}
\ta_{d-1}\left(\Che_{u_d}(\la,z)\right) L^1(z_1)\cdots L^d(z_d)=\\=
L^d(z_d)L^1(z_1)\cdots L^{d-1}(z_{d-1})
\ta_{d-1}\left(\Che_{u_d}(\rh+h^{\le d},z)\right).
\end{multline}
For notational simplicity, set $\la'=\la-h^{>d}$.
\[
\Che_{u_d}(\la,z)=\R^{d-1,d}(\la',\frac{z_{d-1}}{z_d})
\R^{d-2,d}(\la'-h^{d-1},\frac{z_{d-2}}{z_d})\cdots
\R^{1,d}(\la'-h^{\left[ 2,d-1\right]},\frac{z_1}{z_d}),
\]
where $h^{\left[ a,b\right]}$ means $\sum_{a\le j\le b} h^j$. Thus
\begin{equation}\label{eq:tauWphicalc}
\ta_{d-1}\left(\Che_{u_d}(\la,z)\right)=\R^{1,d}(\la',\frac{z_1}{z_d})
\R^{2,d}(\la'-h^1,\frac{z_2}{z_d})\cdots
\R^{d-1,d}(\la'-h^{\le d-2},\frac{z_{d-1}}{z_d}).
\end{equation}
Using \eqref{eq:RLLR1} and the RLL relation \eqref{eq:RLLwithmodifiedR} 
repeatedly, we obtain \eqref{eq:CLLclaim}.
Now the proposition follows by induction on $d$, using that
\[\Che_{t_{d-1}}(\la,z)L^d(z_d)=L^d(z_d)\Che_{t_{d-1}}(\la+h^d,z)\]
which follows from \eqref{eq:RLLR1}.

\end{proof}

The operator $\Cher(\la,z) := \Che_{t_n}(\la,z)$ is called the \emph{Cherednik operator}.
For an operator $F(z)\in \mathcal{E}_n$ we define its matrix elements
$F(z)_{x_1,\ldots,x_n}^{a_1,\ldots,a_n}\in \Fscr$ by 
\[
F(z)(e_{a_1}\otimes\cdots\otimes e_{a_n}\otimes 1) = 
 \sum_{x_1,\ldots,x_n} e_{x_1}\otimes\cdots\otimes e_{x_n}\otimes 
F(z)_{x_1,\ldots,x_n}^{a_1,\ldots,a_n}.
\]

\begin{prop} \label{prop:C1ton}
Let 
\begin{equation}\label{eq:alfatilde}
\tilde\al(\la,z)=\lim_{w\to z}\th(q^2w)\al(\la,w)=
\th(z)\th(q^{2(\la+1)})/\th(q^{2\la}).
\end{equation}
Then
\[\Cher(\la,z)_{1,\ldots, n}^{1,\ldots, n}=
\prod_{i<j}\tilde\al(\la_{ij},z_i/z_j)=
\prod_{i<j}q\th(z_i/z_j)\cdot
\frac{F_{\oneton}(\la)}{F^{\oneton}(\la)}.
\]
\end{prop}
\begin{proof}
The second equality follows from the definition, \eqref{eq:FIdef}, of $F_I$ and $F^I$.
We prove by induction on $d$ that
$\Che_{t_d}(\la,z)_{1,\ldots, n}^{1,\ldots, n}=
\prod_{i<j\le d}\tilde\al(\la_{ij},z_i/z_j)$.
For $d=2$ we have $t_d=s_1$ and
$\Che_{s_1}(\la,z)_{1,\ldots,n}^{1,\ldots,n}=\R^{12}(\la-h^{>2},z_1/z_2)_{1,\ldots,n}^{1,\ldots,n}=
\tilde\al(\la_{12},z_1/z_2)$ as claimed. For $d>2$, 
using the factorization \eqref{eq:chefactorization} we have
\begin{equation}\label{eq:tauWfactor2}
\Che_{t_d}(\la,z)_{1,\ldots,n}^{1,\ldots,n}=
\sum_{x_1,\ldots,x_n}\Che_{t_{d-1}}(\la,z)_{1,\ldots,n}^{x_1,\ldots,x_n}
\ta_{d-1}(\Che_{u_d}(\la,z))_{x_1,\ldots,x_n}^{1,\ldots,n}.
\end{equation}

Since $\Che_{t_{d-1}}(\la,z)$ is a product of operators of the form
$\si(\R^{ii+1}(\la,z_i/z_{i+1}))$ where $1\le i\le d-2$ and
$\si\in S_n$, $\si(j)=j, j>d-1$, and each of these operators
preserve the subspace spanned by
$e_{\ta(1)}\otimes\cdots\otimes e_{\ta(d-1)}\otimes e_d\otimes\cdots\otimes e_n\otimes a$
where $\ta\in S_{d-1}$ and $a\in \Fscr$, the operator $\Che_{t_{d-1}}(\la,z)$ also
preserves this subspace. This means that
$\Che_{t_{d-1}}(\la,z)_{1,\ldots,n}^{x_1,\ldots,x_n}=0$
unless $x_j=j$ for $j\ge d$ and $\{x_1,\ldots,x_{d-1}\}=\{1,\ldots,d-1\}$.
Furthermore, by \eqref{eq:tauWphicalc}, 
\begin{multline}\label{eq:tauWphisum}
\ta_{d-1}(\Che_{u_d}(\la,z))_{x_1,\ldots,x_{d-1},d,\ldots,n}^{1,\ldots,n}=\\=
\sum_{y_2,\ldots, y_{d-1}}\tilde R_{x_1d}^{1y_2}(\la,\frac{z_1}{z_d})
\tilde R_{x_2 y_2}^{2 y_3}(\la-\om(1),\frac{z_2}{z_d})\cdots
\tilde R_{x_{d-1} y_{d-1}}^{d-1,d}(\la-\sum_{k\le d-2}\om(k),\frac{z_{d-1}}{z_d}).
\end{multline}
Here $\tilde R_{xy}^{ab}(\la,z)=\lim_{w\to z} \th(q^2w) R_{xy}^{ab}(\la,w)$.
Since $\tilde R_{xy}^{ab}(\la,z)=0$ unless $\{x,y\}=\{a,b\}$, we deduce that,
when $\{x_1,\ldots,x_{d-1}\}=\{1,\ldots,d-1\}$, 
the terms in the sum \eqref{eq:tauWphisum} are zero unless $x_i=i$ for all $i$
and $y_j=d$ for all $j$. Substituting into \eqref{eq:tauWfactor2}
the claim follows by induction.
\end{proof}

\begin{lem}\label{lem:tdextractable}
 Fix $2 \le d\le n$ and $i<d$. Then there are elements
  $b,c\in B_n$ such that $t_d=s_ib$ and $t_d=cs_i$.
\end{lem}
\begin{proof}
 Since $t_2=s_1$ and $t_3=s_1s_2s_1=s_2s_1s_2$, the statement clearly
 holds for $d=2,3$.
Assuming $d>3$, we first prove the existence of $b$.
 If $i<d-1$ then by induction there is a
  $b'\in B_n$ such that $t_{d-1}=s_ib'$. Hence
  $t_d=t_{d-1}s_{d-1}\cdots s_1=s_i b' s_{d-1}\cdots s_1$.
  Thus we can take $b=b's_{d-1}\cdots s_1$.
If $i=d-1$, write
$t_d=t_{d-2}s_{d-2}\cdots s_{1} s_{d-1}\cdots s_1$.
Then move each of the $d-1$ rightmost factors $s_{d-1},\ldots,s_1$
as far to the left as possible, using that $s_js_k=s_ks_j$ when $|j-k|>1$.
This gives
\[t_d=t_{d-2}s_{d-2}s_{d-1}s_{d-3}s_{d-2}s_{d-4}\cdots s_2s_3s_1s_2s_1.\]
Then use $s_js_{j+1}s_j=s_{j+1}s_js_{j+1}$ repeatedly, working from right
to left, to obtain
\[t_d=t_{d-2}s_{d-1}s_{d-2}s_{d-1}s_{d-3}s_{d-2}\cdots s_4s_2s_3s_1s_2.\]
Finally, $s_{d-1}$ can be moved to the left of $t_{d-2}$ since the latter is a product
of $s_j$'s with $j\le d-3$.

To prove the existence of $c$ we note that $B_n$ carries an involution
$\ast:B_n\to B_n$ satisfying $(a_1a_2)^\ast=a_2^\ast a_1^\ast$ for any $a_1,a_2\in B_n$,
defined by $s_j^\ast=s_j$ for $j\in\oneton$ and $1^\ast=1$. Thus it suffices
to show that $t_d^\ast=t_d$ for any $d$. This is trivial for $d=2,3$.
When $d>3$ we have, by induction on $d$,
\begin{align*}
t_d^\ast&=(t_{d-1}s_{d-1}\cdots s_1)^\ast=s_1\cdots s_{d-1}t_{d-1}=\\
&= s_1\cdots s_{d-1}t_{d-2}s_{d-2}\cdots s_1=\\
&= s_1\cdots s_{d-2}t_{d-2}s_{d-1}s_{d-2}\cdots s_1=\qquad\text{(since $s_{d-1}$ commutes with $t_{d-2}$)}\\
&= t_{d-1}^\ast s_{d-1}\cdots s_1=t_d.
\end{align*}

\end{proof}

\begin{prop}
Let $w=(z_0,q^2z_0,\ldots,q^{2(n-1)}z_0)$, where $z_0\neq 0$ is arbitrary,
 and let $\si,\ta\in S_n$. Then
\begin{equation}\label{eq:cheredniksignidentity}
\Cher(\la,w)_{\si(1),\ldots, \si(n)}^{\ta(1),\ldots, \ta(n)}=
 \frac{\sgn_{\oneton}(\si; \la)}{\sgn^{\oneton}(\ta; \la)}
 \Cher(\la,w)_{1,\ldots, n}^{1,\ldots, n}.
\end{equation}
\end{prop}
\begin{proof}
First we claim that for all $\si,\ta\in S_n$ and each $i\in\oneton$,
\begin{equation}\label{eq:cheredsignpart1}
\Che_{s_i}(\la,w)_{\si\si_i(1),\ldots,\si\si_i(n)}^{\ta(1),\ldots,\ta(n)}=
\si(\sgn_{\oneton}(\si_i;\la))
\Che_{s_i}(\la,w)_{\si(1),\ldots,\si(n)}^{\ta(1),\ldots,\ta(n)},
\end{equation}
and
\begin{equation}\label{eq:cheredsignpart2}
\begin{split}
\Che_{s_i}(\la,w)^{\ta\si_i(1),\ldots,\ta\si_i(n)}_{\si(1),\ldots,\si(n)}&=
\ta(\sgn^{\oneton}(\si_i;\la))
\Che_{s_i}(\la,w)^{\ta(1),\ldots,\ta(n)}_{\si(1),\ldots,\si(n)}=\\
&=-\Che_{s_i}(\la,w)^{\ta(1),\ldots,\ta(n)}_{\si(1),\ldots,\si(n)}.
\end{split}
\end{equation}
Indeed, assume that $z_i/z_{i+1}=q^{-2}$ and that
$\{a_1,\ldots,a_n\}=\{b_1,\ldots,b_n\}=\oneton$.
Then $\Che_{s_i}(\la,z)_{a_1,\ldots,a_n}^{b_1,\ldots,b_n}\neq 0$
iff $\{a_i,a_{i+1}\}=\{b_i,b_{i+1}\}$ in which case
\begin{equation}\label{eq:tildechesiiw}
\Che_{s_i}(\la,z)_{a_1,\ldots,a_n}^{b_1,\ldots,b_n}=
\frac{E(1)E(\la_{a_ia_{i+1}}+1)}{E(\la_{b_{i+1}b_i})}.
\end{equation}
From this and the definitions of the sign functions,
\eqref{eq:leftsigndef}-\eqref{eq:rightsigndef}, the claims follow.
Next, we prove \eqref{eq:cheredniksignidentity} by induction on the sum $\ell$ of 
the lengths of $\si$ and $\ta$. If $\ell=0$ it is trivial. Assuming
\eqref{eq:cheredniksignidentity} holds for $(\si,\ta)$ we prove it holds for
$(\si\si_i,\ta)$ and $(\si,\ta\si_i)$ where $i$ is arbitrary.

Let $i\in\oneton$. By Lemma \ref{lem:tdextractable} we have $t_n=s_ib$ for some $b\in B_n$. 
We have
\begin{align*}
\Che_{t_n}&(\la,w)_{\si\si_i(1),\ldots,\si\si_i(n)}^{\ta(1),\ldots,\ta(n)}=
 \big(\Che_{s_i}(\la,w)\si_i\left(\Che_{b}(\la,w)\right)\big)
 _{\si\si_i(1),\ldots,\si\si_i(n)}^{\ta(1),\ldots,\ta(n)}=\\
&=\sum_{x_1,\ldots,x_n}\Che_{s_i}(\la,w)_{\si\si_i(1),\ldots,\si\si_i(n)}^{x_1,\ldots,x_n}
  \si_i\left(
\Che_{b}(\la,w)\right)_{x_1,\ldots,x_n}^{\ta(1),\ldots,\ta(n)}.
\end{align*}
As in the proof of Proposition \ref{prop:C1ton}, 
$\Che_{s_i}(\la,w)_{\si\si_i(1),\ldots,\si\si_i(n)}^{x_1,\ldots,x_n}$
is zero if $x_1,\ldots,x_n$ is not a permutation of $1,\ldots,n$.
Using \eqref{eq:cheredsignpart1} we obtain
\begin{align*}
&\si(\sgn_{\oneton}(\si_i;\la))
\sum_{x_1,\ldots,x_n}\Che_{s_i}(\la,w)_{\si(1),\ldots,\si(n)}^{x_1,\ldots,x_n}
  \si_i\left(
\Che_{b}(\la,w)\right)_{x_1,\ldots,x_n}^{\ta(1),\ldots,\ta(n)}=\\
&=
\si(\sgn_{\oneton}(\si_i;\la))
\Che_{t_n}(\la,w)_{\si(1),\ldots,\si(n)}^{\ta(1),\ldots,\ta(n)}.
\end{align*}
Using the induction hypothesis and the relation
$\sgn_{\oneton}(\si;\la)
\si(\sgn_{\oneton}(\si_i;\la))=\sgn_{\oneton}(\si\si_i;\la)$
we obtain \eqref{eq:cheredniksignidentity} for $(\si\si_i,\ta)$.

For the other case, let $i$ be arbitrary and set $j=\ta_n(i)$.
By Lemma \ref{lem:tdextractable} there is a $c\in B_n$ such that $t_n=cs_j$.
Recall the surjective morphism $\pi:B_n\to S_n$ sending $s_i$ to $\si_i=(i\; i+1)$.
Then $\si_j\pi(c)=\pi(c)\si_i$. We have
\begin{align*}
\Che_{t_n}&(\la,w)_{\si(1),\ldots,\si(n)}^{\ta\si_i(1),\ldots,\ta\si_i(n)}=
\big(\Che_c(\la,w)\cdot\pi(c)(\Che_{s_j}(\la,w))\big)_{\si(1),\ldots,
 \si(n)}^{\ta\si_i(1),\ldots\ta\si_i(n)}=\\
 &=
\sum_{x_1,\ldots,x_n}
\Che_c(\la,w)_{\si(1),\ldots,\si(n)}^{x_1,\ldots,x_n}
  \pi(c)\left(\Che_{s_j}(\la,w)\right)_{x_1,\ldots,x_n}^{\ta\si_i(1),\ldots,\ta\si_i(n)}.
\end{align*}
It is easy to check that $\si(F(z))_{a_1,\ldots,a_n}^{b_1,\ldots,b_n}=
F(\si(z))_{a_{\si(1)},\ldots,a_{\si(n)}}^{b_{\si(1)},\ldots,b_{\si(n)}}$
for any $F(z)\in\mathcal{E}_n$ and $\si\in S_n$.
Define $w_i$ by $(w_1,\ldots,w_n)=w=(z_0,q^2z_0,\ldots,q^{2(n-1)}z_0)$.
Then $w_i/w_{i+1}=q^{-2}$ for each $i$.
Set $w'=(w_{\pi(c)(1)},\ldots,w_{\pi(c)(n)})$.
For each $i$,
$w_{\pi(c)(i)}/w_{\pi(c)(i+1)}=w_{\ta_n(i+1)}/w_{\ta_n(i)}=q^{-2}$ also.
Therefore
\begin{align*}
\Che_{t_n}&(\la,w)_{\si(1),\ldots,\si(n)}^{\ta\si_i(1),\ldots,\ta\si_i(n)}=
\sum_{x_1,\ldots,x_n} \Che_c(\la,w)_{\si(1),\ldots,\si(n)}^{x_1,\ldots,x_n}
\Che_{s_j}(\la,w')_{x_{\pi(c)(1)},\ldots,x_{\pi(c)(n)}}^{\ta\si_i\pi(c)(1),\ldots,
\ta\si_i\pi(c)(n)}=\\
&=
\sum_{x_1,\ldots,x_n} \Che_c(\la,w)_{\si(1),\ldots,\si(n)}^{x_1,\ldots,x_n}
\Che_{s_j}(\la,w')_{x_{\pi(c)(1)},\ldots,x_{\pi(c)(n)}}^{\ta\pi(c)\si_j(1),\ldots,
\ta\pi(c)\si_j(n)}=\\
&=
\sum_{x_1,\ldots,x_n} \Che_c(\la,w)_{\si(1),\ldots,\si(n)}^{x_1,\ldots,x_n}
(\sgn \si_j) \Che_{s_j}(\la,w')_{x_{\pi(c)(1)},\ldots,x_{\pi(c)(n)}}^{\ta\pi(c)(1),\ldots,
\ta\pi(c)(n)}=\\
&=
\sum_{x_1,\ldots,x_n} \Che_c(\la,w)_{\si(1),\ldots,\si(n)}^{x_1,\ldots,x_n}
(-1) \pi(c)(\Che_{s_j}(\la,w)_{x_1,\ldots,x_n}^{\ta(1),\ldots,\ta(n)}=\\
&=
- \Che_{t_n}(\la,w)_{\si(1),\ldots,\si(n)}^{\ta(1),\ldots,\ta(n)}.
\end{align*}
By the induction hypothesis it follows that
\eqref{eq:cheredniksignidentity} holds for $(\si,\ta\si_i)$.
This proves the formula \eqref{eq:cheredniksignidentity}.
\end{proof}

\begin{thm} \label{thm:leftequalsright}
For any subsets $I,J\subseteq \oneton$ and $z\in\C^\times$, the left and
right elliptic minors coincide:
\[\overleftarrow{\xi}_I^J(z)=\overrightarrow{\xi}_I^J(z).\]
We denote this common element by $\xi_I^J(z)$.
\end{thm}
\begin{proof}
If $\#I\neq \#J$ both sides are zero. Suppose $\#I=\#J=d$.
By relation \eqref{eq:Snonminors} we can, after applying a suitable automorphism,
assume that $I=J=\onetod$.
Since the subalgebra of $\Fscr$ generated by $e_{ij}(z)$, $i,j\in\onetod$, $z\in\C^\times$
and $f(\la), f(\rh)$ with $f\in M_{\hf_d^*}\subseteq M_{\hf^*}$, $\hf_d$ being the Cartan
subalgebra of $\mathfrak{sl}(d)$,
is isomorphic to $\Fell(M(d))$, we can also assume $d=n$.
Identifying the matrix element $_{1,\ldots,n}^{1,\ldots,n}$
on both sides of \eqref{eq:WLLrel} we get
\begin{multline*}
\sum_{x_1,\ldots,x_n} \Cher(\la,z)_{1,\ldots,n}^{x_1,\ldots,x_n}
 e_{x_1, 1}(z_1)\cdots e_{x_n, n}(z_n)=\\=
\sum_{x_1,\ldots,x_n} 
 e_{n,x_n}(z_n)\cdots e_{1,x_1}(z_1)
 \Cher(\rh+h^{\le n},z)_{x_1,\ldots,x_n}^{1,\ldots,n}.
\end{multline*}
As in the proof of Proposition \ref{prop:C1ton},
$\Cher(\la,z)_{1,\ldots,n}^{x_1,\ldots,x_n}$ is zero if $x_1,\ldots,x_n$ is not a
permutation of $1,\ldots,n$.
Taking $z=w=(z_0,q^2z_0,\ldots,q^{2(n-1)}z_0)$ and dividing both sides by
$\prod_{i<j}q\th(w_i/w_j)=\prod_{i<j}q\th(q^{2(i-j)})$ we get
\begin{multline*}
\frac{F_{\oneton}(\la)}{F^{\oneton}(\la)}
\sum_{\si\in S_n} \sgn^{\oneton}(\si;\la)^{-1}
 e_{\si(1)1}(z_0)\cdots e_{\si(n)n}(q^{2(n-1)}z_0)=\\
=\frac{F_{\oneton}(\rh)}{F^{\oneton}(\rh)}
\sum_{\ta\in S_n} 
 \sgn_{\oneton}(\ta;\rh)
 e_{n\ta(n)}(q^{2(n-1)}z_0)\cdots e_{1\si(1)}(z_0)
\end{multline*}
Multiplying by $\frac{F^{\oneton}(\rh)}{F_{\oneton}(\la)}$ and comparing
with \eqref{eq:leftminor} and \eqref{eq:rightminor}, we deduce that
$\overrightarrow{\xi}_{\oneton}^{\oneton}(z_0)=\overleftarrow{\xi}_{\oneton}^{\oneton}(z_0)$,
as desired.
\end{proof}

\subsection{Laplace expansions}\label{sec:laplace}
Using the left (right) $\Fell(M(n))$-comodule algebra structure of $\La$
($\La'$) it is straightforward to prove Laplace expansion formulas for the
elliptic minors.
For subsets $I,J\subseteq\oneton$ we define $S_l(I,J;\ze), S_r(I,J;\ze)
\in M_{\hf^*}$ by
\begin{align}
 \label{eq:Sldef}
 v_I(q^{2\# J}z)v_J(z)&=S_l(I,J;\ze)v_{I\cup J}(z),\\
 \label{eq:Srdef}
 w^I(z)w^J(q^{2\# I}z)&=S_r(I,J;\ze)w^{I\cup J}(z).
\end{align}
That this is possible follows from the definitions \eqref{eq:vIdef},\eqref{eq:wIdef}
of $v_I(z), w^I(z)$ and the commutation relations \eqref{eq:leftsimp1}-\eqref{eq:leftsimp2},
\eqref{eq:rightsimp1}-\eqref{eq:rightsimp2}.
In particular $S_l(I,J;\ze)=0=S_r(I,J;\ze)$ if $I\cap J\neq\emptyset$.
\begin{thm}
(i) Let $I_1,I_2,J\subseteq\oneton$ and set $I=I_1\cup I_2$. Then
\begin{equation}\label{eq:leftLaplace}
 S_l(I_1,I_2;\la)\xi_I^J(z)=\sum_{J_1 \cup J_2 = J} S_l(J_1,J_2;\rh)
 \xi_{I_1}^{J_1}(q^{2\# I_2}z)\xi_{I_2}^{J_2}(z).
\end{equation}
(ii) Let $J_1,J_2,I\subseteq\oneton$ and set $J=J_1\cup J_2$. Then
\begin{equation}\label{eq:rightLaplace}
S_r(J_1,J_2;\rh)\xi_I^J(z)=\sum_{I_1\cup I_2=I} S_r(I_1,I_2;\la)
\xi_{I_1}^{J_1}(z)\xi_{I_2}^{J_2}(q^{2\# J_1}z).
\end{equation}
\end{thm}
\begin{proof}
 We have
 \begin{multline*}
 \De_{\La}(v_{I_1}(q^{2\# I_2}z))\De_{\La}(v_{I_2}(z))=
 \sum_{J_1,J_2} \xi_{I_1}^{J_1}(q^{2\# I_2}z)\xi_{I_2}^{J_2}(z)\otimes 
 v_{J_1}(q^{2\# I_2}z) v_{J_2}(z)=\\
 =\sum_{J_1,J_2} \xi_{I_1}^{J_1}(q^{2\# I_2}z)\xi_{I_2}^{J_2}(z)\otimes
 S_l(J_1,J_2;\ze)v_J(z)=\\
 =\sum_J \left(\sum_{J_1\cup J_2=J} S_l(J_1,J_2;\rh)\xi_{I_1}^{J_1}(q^{2\# I_2}z)
 \xi_{I_2}^{J_2}(z)\right) \otimes v_J(z).
\end{multline*}
On the other hand,
\begin{multline*}
 \De_{\La}(v_{I_1}(q^{2\# I_2}z))\De_{\La}(v_{I_2}(z))=
 \De_{\La}(v_{I_1}(q^{2\# I_2}z)v_{I_2}(z))=\\
 =\De_{\La}(S_l(I_1,I_2;\ze)v_I(z))=\sum_J S_l(I_1,I_2;\la)\xi_I^J(z)\otimes v_J(z).
\end{multline*}
Equating these expressions proves \eqref{eq:leftLaplace} since,
by Proposition \ref{prop:extalg},
the set $\{v_J(z):J\subseteq\oneton\}$
is linearly independent over $M_{\hf^*}$.
The second part is completely
analogous, using the right comodule algebra $\La'$ in place of $\La$.
\end{proof}

In Section \ref{sec:antipode} we will need the following lemma,
relating the left and right signums
$S_l(I,J;\ze)$ and $S_r(I,J;\ze)$, defined in \eqref{eq:Sldef},\eqref{eq:Srdef}.
In the non-spectral trigonometric case
the corresponding identity was proved in \cite{N05}
(in proof of Proposition 4.1.22).

\begin{lem} Let $I,J$ be two disjoint subsets of $\oneton$. Then
\begin{equation}\label{eq:SlSrrelation}
 S_l(I,J;\ze+\om(I))=S_r(J,I;\ze)^{-1}
\end{equation}
where $\om(I)=\sum_{i\in I}\om(i)$.
\end{lem}
\begin{proof}
First we claim that, we have the following explicit formulas:
\begin{align}
\label{eq:Slformula}
S_l(I,J;\ze)&=\prod_{i\in I,j\in J} E(\ze_{ji}+1),\\
\label{eq:Srformula}
S_r(I,J;\ze)&= \prod_{i\in I,j\in J} E(\ze_{ij})^{-1}.
\end{align}
Recall the definition, \eqref{eq:vIdef}, of $v_I(z)$.
Since $E$ is odd, relation \eqref{eq:leftsimp1} implies that
\[v_i(q^2z)v_j(z)=\frac{E(\ze_{ji}+1)}{E(\ze_{ij}+1)}v_j(q^2z)v_i(z).\]
Also, $F_J(\ze)$ only involves $\ze_{ij}$ with $i,j\in J$ so
it commutes with any $v_k(z)$ with $k\in I$ (since $I\cap J=\emptyset$).
From these facts we obtain
\begin{align*} 
v_I(q^{2\# J}z)v_J(z)&=\frac{F_I(\ze)^{-1}F_J(\ze)^{-1}}{F_{I\cup J}(\ze)^{-1}}
\prod_{\substack{i\in I,j\in J\\ i<j}}\frac{E(\ze_{ji}+1)}{E(\ze_{ij}+1)}v_{I\cup J}(z)=\\
&=\prod_{\substack{(i,j)\in K\\i<j}} E(\ze_{ij}+1)
\prod_{\substack{i\in I,j\in J\\ i<j}}\frac{E(\ze_{ji}+1)}{E(\ze_{ij}+1)}v_{I\cup J}(z)=\\
&=\prod_{i\in I, j\in J} E(\ze_{ji}+1) v_{I\cup J}(z),
\end{align*}
where $K=(I\times J)\cup (J\times I)$. This proves \eqref{eq:Slformula}. Similarly
one proves \eqref{eq:Srformula}. Now we have
\[S_l(J,I;\ze+\om(J))^{-1}=\prod_{i\in I,j\in J} E((\ze+\om(J))_{ij}+1)^{-1}=
\prod_{i\in I,j\in J}E(\ze_{ij})^{-1}=S_r(I,J;\ze).\]
Here we used that for any $i\in I, j\in J$ we have $\om(J)(E_{ii})=0$ and
$\om(J)(E_{jj})=1$ and hence $(\om(J))_{ij}=-1$.
\end{proof}

\section{The cobraiding and the elliptic determinant}\label{sec:pairingsanddet}

\subsection{Cobraidings for $\hf$-bialgebroids} \label{sec:pairings}
The following definition of a cobraiding was given in \cite{Ro}.
When $\hf=0$ the notion reduces to ordinary cobraidings for bialgebras.
\begin{dfn}
 A \emph{cobraiding} on an $\hf$-bialgebroid $A$ is a $\C$-bilinear map
 $\langle\cdot,\cdot\rangle:A\times A\to D_\hf$ such that,
 for any $a,b,c\in A$ and $f\in M_{\hf^\ast}$,
\begin{subequations}
\begin{align}
\label{eq:cobraidinghomogen}
 &\langle A_{\al\be}, A_{\ga\de}\rangle\subseteq (D_\hf)_{\al+\ga,\be+\de},\\
 \label{eq:pairing_mu1}
 &\langle \mu_r(f)a,b\rangle = \langle a,\mu_l(f)b\rangle = fT_0 \circ \langle a,b\rangle,\\
 \label{eq:pairing_mu2}
 &\langle a\mu_l(f),b\rangle = \langle a,b\mu_r(f)\rangle = \langle a,b\rangle \circ fT_0,\\
 \label{eq:pairing_ab_c}
 &\langle ab,c\rangle = \sum_i \langle a,c_i'\rangle T_{\be_i}\langle b,c_i''\rangle,
 \quad \De(c)=\sum_i c_i'\otimes c_i'',\quad c_i''\in A_{\be_i\ga},\\
 \label{eq:pairing_a_bc}
 &\langle a,bc\rangle = \sum_i \langle a_i'',b\rangle T_{\be_i} \langle a_i',c\rangle,
 \quad \De(a)=\sum_i a_i'\otimes a_i'',\quad a_i''\in A_{\be_i\ga},\\
 &\langle a,1\rangle = \langle 1,a\rangle = \ep(a),\\
 &\sum_{ij} \mu_l\big( \langle a_i', b_j'\rangle 1\big) a_i'' b_j'' =
 \sum_{ij} \mu_r\big( \langle a_i'', b_j ''\rangle 1\big) b_j' a_i'.
 \label{eq:cobraidingidentity}
\end{align}
\end{subequations}
\end{dfn}

The following definition was given in unpublished notes by
Rosengren \cite{Rnotes}. The terminology is motivated by Proposition
\ref{prop:unitarycob} below concerning FRST algebras $A_R$, but it
makes sense for arbitrary $\hf$-bialgebroids.
\begin{dfn}
A cobraiding $\langle\cdot,\cdot\rangle$ on an $\hf$-bialgebroid
$A$ is called \emph{unitary} if
\begin{equation}
\ep(ab)=\sum_{(a),(b)}\langle a',b'\rangle
T_{\om_{12}(a)+\om_{12}(b)}\langle a'',b''\rangle
\end{equation}
for all $a,b\in A$.
\end{dfn}

\subsection{Cobraidings for the FRST-algebras $A_R$.}
Now let $R:\hf^*\times \C^\times \to \End_\hf(V\otimes V)$ be a meromorphic
function and let $A_R$ be the $\hf$-bialgebroid associated to $R$ as in
Section \ref{sec:FRST}.

\begin{prop}\label{prop:pairing}
 Assume that $\ph:\C^\times\to\C$ is a holomorphic function,
 not vanishing identically,
  such that, for each $x,y,a,b\in X$, $z\in\C^\times$,
 the limit $\lim_{w\to z}\left(\ph(w)R_{xy}^{ab}(\ze,w)\right)$ exists and defines a meromorphic
 function in $M_{\hf^*}$. Then the following statements are equivalent:
 \begin{itemize}
  \item[(i)] there exists a cobraiding $\langle\cdot,\cdot\rangle:A_R\times A_R\to D_\hf$
  satisfying
  \begin{equation}\label{eq:pairingFRSTgammal}
   \langle L_{ij}(z_1),L_{kl}(z_2)\rangle = \lim_{w\to z_1/z_2}\left(\ph(w)R_{ik}^{jl}(\ze,w)\right)
   T_{-\om(i)-\om(k)},
  \end{equation}
  \item[(ii)] $R$ satisfies the QDYBE \eqref{eq:qdybe}.
 \end{itemize}
\end{prop}
\begin{rem}
a) The identity \eqref{eq:cobraidingidentity} is not necessary when proving that (i) implies (ii).
Without assuming \eqref{eq:cobraidingidentity}, $\langle \cdot,\cdot\rangle$ is a
\emph{pairing on $A^{\text{cop}}\times A$}. See \cite{Ro}.

b)
Without the factor $\ph(w)$, the cobraiding is
not well-defined if $R(\ze,z)$ has poles in the $z$ variable.
We also remark that the residual relations \eqref{eq:RLLmatrixelements_residual_general}
are necessary for (ii) to imply (i).
\end{rem}
\begin{proof}
The proof is straightforward and is carried out in \cite{N05}, Lemma 2.2.5,
under the assumption that the R-matrix is regular in the spectral variable.
\end{proof}

We will now generalize slightly the notion of a unitary cobraiding on
$A_R$ to account for spectral singularities in the R-matrix as follows.

Call $a\in A_R$ \emph{spectrally homogenous} if there
exist $k\in\Z_{\ge 0}$ and $z_1,\ldots,z_k\in\C^\times$ such that
\begin{equation}
a\in\sum_{\si\in S_k}\sum_{i_l, j_l\in X} M_{\hf^\ast}\otimes M_{\hf^\ast}
L_{i_1j_1}(z_{\si(1)})\cdots L_{i_kj_k}(z_{\si(k)}).
\end{equation}
The multiset $\{z_i\}_i$ is called the \emph{spectral degree} of $a$
and is denoted $\sdeg(a)$.
Note that the spectral degree of a nonzero spectrally homogenous
element is uniquely defined, since the RLL-relations 
\eqref{eq:RLLmatrixelements_residual_general} are spectrally homogenous.

Let $\ph:\C^\times\to\C$ be holomorphic.
For spectrally homogenous elements $a,b\in A_R$, define the
\emph{regularizing factor} $\widehat{\ph}(a,b)$ by
\begin{equation}
\widehat{\ph}(a,b)=\prod_{1\le i\le k, 1\le j\le l} \ph(z_i/w_j)
\end{equation}
where $\{z_i\}_i=\sdeg(a)$, $\{w_j\}_j=\sdeg(b)$.

\begin{dfn}
Let $\ph:\C^\times\to\C$ be holomorphic.
A cobraiding $\langle\cdot,\cdot\rangle$ on $A_R$ is \emph{unitary with
respect to $\ph$} if
\begin{equation}\label{eq:unitarycobwrtph}
\widehat{\ph}(a,b)\widehat{\ph}(b,a)\ep(ab)=\sum_{(a),(b)}\langle a',b'\rangle
T_{\om_{12}(a)+\om_{12}(b)}\langle a'',b''\rangle
\end{equation}
for all spectrally homogenous $a,b\in A_R$.
\end{dfn}

The following proposition was proved in \cite{Rnotes}
if the spectral variables is taken to be generic
so that no regularizing factors are needed.
\begin{prop}\label{prop:unitarycob}
Suppose $R:\hf^\ast\times\C^\times\to\End_\C(V\otimes V)$ satisfies
the QDYBE and is unitary: $R(\ze,z)R(\ze,z^{-1})^{(21)}=\Id_{V\otimes V}$.
Suppose $\ph:\C^\times\to\C$ is nonzero holomorphic such that
$\lim_{w\to z}\left(\ph(w)R_{xy}^{ab}(\ze,w)\right)$ exists and
is a holomorphic function in $M_{\hf^\ast}$.
Then the cobraiding $\langle\cdot,\cdot\rangle$ on $A_R$ given in Proposition
\ref{prop:pairing} is unitary with respect to $\ph$.
\end{prop}
\begin{proof}
Since both sides are holomorphic in the spectral variables,
it is enough to prove it for generic values. We claim that for such values,
$\widehat{\ph}(a,b)\langle a,b\rangle_R = \langle a,b\rangle$
where $\langle\cdot,\cdot\rangle_R$ is the cobraiding, defined only
for generic spectral values, determined by
$\langle L_{ij}(z),L_{kl}(w)\rangle_R = R_{ik}^{jl}(\ze,z/w)T_{-\om(i)-\om(k)}$.
Indeed, this claim follows by induction from the identities
\eqref{eq:pairing_ab_c},\eqref{eq:pairing_a_bc} using that
$\widehat{\ph}(a_1,b)\widehat{\ph}(a_2,c)=\widehat{\ph}(a_3,bc)$ and
$\widehat{\ph}(c,a_1)\widehat{\ph}(b,a_2)=\widehat{\ph}(cb,a_3)$
for spectrally homogenous $a_i,b,c\in A_R$, the $a_i$ having the
same spectral degree.

Since the R-matrix $R$ is unitary, the statement of the lemma now
follows from the identity
\begin{equation}
\ep(ab)=\sum_{(a),(b)}\langle a',b'\rangle_R
T_{\om_{12}(a)+\om_{12}(b)}\langle a'',b''\rangle_R
\end{equation}
holding for generic spectral values which was proved by Rosengren \cite{Rnotes}.
\end{proof}

\subsubsection{The case of $\Fell(M(n))$}
Specializing further to the algebra of interest, $\Fell(M(n))$, we obtain
the following corollary.
\begin{cor}
The $\hf$-bialgebroid $\Fell(M(n))$ carries a cobraiding $\langle\cdot,\cdot\rangle$
satisfying
\begin{equation}
\langle e_{ij}(z), e_{kl}(w) \rangle = \tilde R_{ik}^{jl}(\zeta,z/w)
 T_{-\om(i)-\om(k)}\qquad\forall z,w\in\C^\times, i,j\in[1,n],
\end{equation}
where
\begin{equation}\label{eq:cobraidingFRST}
\tilde R_{ik}^{jl}(\zeta,z) = \lim_{w\to z} \big(\th(q^2 w)R_{ik}^{jl}(\zeta,w)\big).
\end{equation}
Moreover, this cobraiding is unitary with respect to
$\ph:\C^\times\to\C$, $\ph(z)=\th(q^2z)$.
\end{cor}
\begin{proof}
It suffices to notice that, by
\eqref{eq:Rmatrixelts}, \eqref{eq:aldef},\eqref{eq:bedef},
 $\tilde R$ is regular in $z$, and apply Proposition \ref{prop:pairing}
 and Proposition \ref{prop:unitarycob}.
\end{proof}

\subsection{Properties of the elliptic determinant} \label{sec:centrality}

A common method used to study quantum minors and prove
that quantum determinants are central
is the fusion procedure, going back to work by
Kulish and Sklyanin \cite{KS82}. Another approach, using
representation theory, was developed by Noumi, Yamada and Mimachi
\cite{NYM}.
In this section we show how to prove that the elliptic determinant is central
using the properties of the cobraiding on $\Fell(M(n))$ and how to resolve
technical issues connected with the spectral singularities of
the elliptic R-matrix.

Let $A=\Fell(M(n))$. When $I=J=\oneton$ we set 
\begin{equation}
 \det(z)=\xi_I^J(z)
\end{equation}
for $z\in\C^\times$, where $\xi_I^J(z)$ is the elliptic minor given in
Theorem \ref{thm:leftequalsright}. Thus one possible expression for $\det(z)$ is
\begin{equation}
\det(z)=\sum_{\si\in S_n}
 \frac{ F^{\oneton}(\rh) }{ \si\big(F^{\oneton}(\la)\big)}
 e_{\si(1)1}(z)e_{\si(2)2}(q^2z)\cdots e_{\si(n)n}(q^{2(n-1)}z).
\end{equation}

\begin{thm}\label{thm:det}
a) $\det(z)$ is a grouplike element of $A$ for each $z\in\C^\times$, i.e.
\[\De(\det(z))=\det(z)\otimes\det(z),\qquad \ep(\det(z))=1.\]
b) $\det(z)$ is a central element in $\Fell(M(n))$:
\begin{equation}
[e_{ij}(z),\det(w)]=[f(\la),\det(w)]=[f(\rh),\det(w)]=0
\end{equation}
for all $f\in M_{\hf^\ast}$, $i,j\in [1,n]$ and all $z,w\in\C^\times$.
\end{thm}
\begin{proof}
Let $\La^n(z)=M_{\hf^*} v_I(z)$, where $I=\oneton$. It is a one-dimensional
subcorepresentation of the left exterior corepresentation $\La$. Its matrix element
is $\det(z)$, i.e.
\[\De(v_I(z))=\det(z)\otimes v_I(z).\]
From the coassociativiy and counity axioms for a corepresentation
follows that $\det(z)$ is grouplike, proving part a).

The rest of this section is devoted to the proof of part b).
It follows from the definition that 
$\det(z)\in A_{00}$ and thus it commutes with $f(\rho)$
and $f(\la)$ for any $f\in M_{\hf^*}$.
To prove that it commutes with the generators $e_{ij}(z)$ we need
several lemmas which we now state and prove.

\begin{lem} \label{lem:pairingzeros}
For $i,j\in[1,n]$, $I,J\subseteq[1,n]$, $\#I=\#J=2$, we have
\begin{alignat}{2}
\label{eq:pairingz=w}
\left\langle \xi_I^J(w),e_{ij}(z) \right\rangle = 0 &\quad\text{if $w/z\in p^\Z$,}\\
\label{eq:pairingz=q^2w}
\left\langle e_{ij}(z),\xi_I^J(w) \right\rangle = 0 &\quad\text{if $q^2w/z\in p^\Z$.}
\end{alignat}
\end{lem}
\begin{proof}
Let $I=\{i_1,i_2\}$, $i_1<i_2$, $J=\{j_1,j_2\}$, $j_1<j_2$. Using the
left expansion formula \eqref{eq:leftminor} and
\eqref{eq:pairing_mu1},\eqref{eq:pairing_mu2} we have
\begin{align*}
\left\langle \xi_I^J(w),e_{ij}(z)\right\rangle &=
\left\langle\frac{E(\rho_{j_1j_2}+1)}{E(\la_{i_1i_2}+1)}
e_{i_2j_2}(q^2w)e_{i_1j_1}(w), e_{ij}(z)\right\rangle 
 \quad  +[j_1 \leftrightarrow j_2]=\\
&=E(\zeta_{j_1j_2+1}+1)\left\langle e_{i_2j_2}(q^2w)e_{i_1j_1}(w),e_{ij}(z)
\right\rangle
\frac{1}{E(\zeta_{i_1i_2}+1)} \quad +[j_1\leftrightarrow j_2].
\end{align*}
Thus we need to prove that for $w/z\in p^\Z$,
the first term is anti-symmetric in $j_1,j_2$.
By \eqref{eq:pairing_ab_c},
\begin{align}
&E(\ze_{j_1j_2}+1)\left\langle e_{i_2j_2}(q^2w)e_{i_1j_1}(w), e_{ij}(z)
 \right\rangle =\nonumber\\
&=E(\ze_{j_1j_2}+1)\sum_x\left\langle e_{i_2j_2}(q^2w), e_{ix}(z)\right\rangle
 T_{\om(x)}\left\langle e_{i_1j_1}(w), e_{xj}(z)\right\rangle=\nonumber\\
&=E(\ze_{j_1j_2}+1)\sum_x
 \tilde R_{i_2 i}^{j_2 x}\big(\ze,q^2w/z\big)
 \tilde R_{i_1x}^{j_1j}\big(\ze-\om(j_2),w/z\big)
 T_{-\om(j_1)-\om(j_2)-\om(j)}. \label{eq:pz_expr1}
\end{align}
Take now $w=p^kz$ where $k\in\Z$.
One checks that
\begin{equation}\label{eq:pz_id0}
\tilde R^{ab}_{xy}(\ze,p^k)=\th(q^2p^k)q^{2k(\ze_{ba}+1-\delta_{ab})}\delta_{ay}\delta_{bx},
\end{equation}
where $\delta_{xy}$ is the Kronecker delta.
In particular, only the $x=j_1$ term is nonzero. Now the anti-symmetry
of \eqref{eq:pz_expr1} in $j_1,j_2$ 
follows by applying the identities
\begin{equation}\label{eq:pz_id1}
\tilde R^{j_1j}_{jj_1}\big(\ze-\om(j_2),p^k\big)=
 q^{2k\ze_{j_2j_1}}\cdot
 \tilde R^{j_2j}_{jj_2}\big(\ze-\om(j_1),p^k\big),
\end{equation}
\begin{equation}\label{eq:pz_id2}
E(\ze_{j_1j_2}+1) \tilde R_{i_2 i}^{j_2 j_1}\big(\ze,q^2p^k\big)=
-q^{2k\ze_{j_1j_2}}\cdot
E(\ze_{j_2j_1}+1) \tilde R_{i_2 i}^{j_1 j_2}\big(\ze,q^2p^k\big).
\end{equation}
Relation \eqref{eq:pz_id1} can be proved directly from \eqref{eq:pz_id0}
while for \eqref{eq:pz_id2} one can use that
\[\frac{\tilde R_{xy}^{ab}(\ze,p^kz)}
      {\tilde R_{xy}^{ba}(\ze,p^kz)}=q^{2k\ze_{ba}}
 \frac{\tilde R_{xy}^{ab}(\ze,z)}
      {\tilde R_{xy}^{ba}(\ze,z)}\]
together with the relation    
\[E(\ze_{j_1j_2}+1)\tilde\al(\ze_{j_2j_1},q^2)=-E(\ze_{j_2j_1}+1)
\tilde\be(\ze_{j_2j_1},q^2)\]
which holds for any $j_1\neq j_2$ which is easily proved by
applying $\theta(z^{-1})=-z^{-1}\theta(z)$ three times.

Relation \eqref{eq:pairingz=q^2w} can be proved analogously, using the right 
expansion formula \eqref{eq:rightminor} for $\xi_I^J(w)$ instead.
\end{proof}

Since the cobraiding depend holomorphically on the spectral variables
and all zeros of $\theta$ are simple and of the form $p^k$,
we conclude that the following limits exist for all $z,w\in\C^\times$,
$i,j,I,J$, $\#I=\#J=2$:
\begin{align}
\left\langle \xi_I^J(w), e_{ij}(z)\right\rangle' &:=
 \lim_{(z_1,w_1)\to (z,w)}\frac{\left\langle \xi_I^J(w_1),
  e_{ij}(z_1)\right\rangle}{\theta(z_1/w_1)},\\
\left\langle e_{ij}(z),\xi_I^J(w)\right\rangle' &:=
 \lim_{(z_1,w_1)\to (z,w)}\frac{\left\langle e_{ij}(z_1),
  \xi_I^J(w_1)\right\rangle}{\theta(q^2w_1/z_1)}.
\end{align}

Taking $a=e_{ij}(z_1), b=\xi_I^J(w_1)$
in \eqref{eq:unitarycobwrtph},
dividing both sides by $\theta(z_1/w_1)\theta(q^2w_1/z_1)$ and taking the
limits $(z_1,w_1)\to(z,w)$, where $z,w\in\C^\times$ are arbitrary, we get
\begin{align}
\label{eq:unitarycob3}
\psi(z,w)
 \ep\left(e_{ij}(z)\xi_I^J(w)\right) &=
 \sum_{x,X}\left\langle\xi_I^X(w),e_{ix}(z)\right\rangle'
 T_{\om(x)+\om(x_1)+\om(x_2)}\left\langle e_{xj}(z),\xi_X^J(w)\right\rangle',
\\
\intertext{and interchanging $a$ and $b$,}
\label{eq:unitarycob4}
\psi(z,w)\ep\left(\xi_I^J(w)e_{ij}(z)\right)&=
 \sum_{x,X}\left\langle e_{ix}(z),\xi_I^X(w)\right\rangle'
 T_{\om(x)+\om(x_1)+\om(x_2)}\left\langle\xi_X^J(w),e_{xj}(z)\right\rangle',
\end{align}
for all $z,w\in\C^\times$,
where $\psi:(\C^\times)^2\to\C$ is given by
\begin{equation}
\psi(z,w)=\theta(q^2z/w)
\theta(q^4w/z).
\end{equation}

We are now ready to prove the key identity.
\begin{lem}\label{lem:keyidentity}
For any $i,j\in[1,n]$, $I,J\subseteq [1,n]$, $\#I=\#J=2$ and any
$z,w\in\C^\times$, $q^2w/z\notin p^\Z$ we have
\begin{multline}\label{eq:keyidentity}
\psi(z,w)\sum_{x,X}\mu_l\left(\left\langle\xi_I^X(w),e_{ix}(z)\right\rangle'1\right)
 \xi_X^J(w)e_{xj}(z)=\\=
\psi(z,w)\sum_{x,X}\mu_r\left(\left\langle\xi_X^J(w),e_{xj}(z)\right\rangle'1\right)
 e_{ix}(z)\xi_I^X(w).
\end{multline}
\end{lem}

\begin{proof}
Using the counit axiom followed by \eqref{eq:unitarycob3} we have
\begin{align}
\psi(z,w) e_{ij}(z)\xi_I^J(w)&=
 \psi(z,w)\sum_{x,X}\mu_l\left(\ep\left(e_{ix}(z)\xi_I^X(w)\right)1\right)
 e_{xj}(z)\xi_X^J(w)=\nonumber \\
&=\sum_{x,y,X,Y}\mu_l\left(\left\langle\xi_I^Y(w),e_{iy}(z)\right\rangle' 1\right)
 \mu_l\left(\left\langle e_{yx}(z),\xi_Y^X(w)\right\rangle' 1\right)
 e_{xj}(z)\xi_X^J(w). \label{eq:keylemmaproof1}
\end{align}
Applying the identity obtained
by dividing by $\theta(q^2w/z)$ in both sides of the cobraiding identity
\eqref{eq:cobraidingidentity} with $a=e_{yj}(z), b=\xi_Y^J(w)$ in
the right hand side of \eqref{eq:keylemmaproof1} gives
\[
\psi(z,w) e_{ij}(z)\xi_I^J(w)=
\sum_{x,y,X,Y}
 \mu_l\left(\left\langle \xi_I^Y(w),e_{iy}(z) \right\rangle' 1\right)
 \mu_r\left(\left\langle e_{xj}(z),\xi_X^J(w) \right\rangle' 1\right)
 \xi_Y^X(w)e_{yx}(z).
\]
Now multiply both sides by $\mu_r\left(\left\langle \xi_J^K(w),e_{jk}(z)
\right\rangle' 1\right)$ and sum over $j,J$.
After applying \eqref{eq:unitarycob4} in the right hand side we get
\begin{multline*}
\psi(z,w)\sum_{j,J}\mu_r\left(\left\langle\xi_J^K(w),e_{jk}(z)\right\rangle' 1
\right) e_{ij}(z)\xi_I^J(w)=\\=
\psi(z,w)
\sum_{x,y,X,Y}
 \mu_l\left(\left\langle \xi_I^Y(w),e_{iy}(z) \right\rangle' 1\right)
 \mu_r\left(\ep\left(\xi_X^K(w)e_{xk}(z)\right)1\right)
 \xi_Y^X(w)e_{yx}(z).
\end{multline*}
By the counit axiom the last expression equals
\[\psi(z,w)\sum_{y,Y}
 \mu_l\left(\left\langle \xi_I^Y(w), e_{iy}(z)\right\rangle '1\right)
 \xi_Y^K(w)e_{yk}(z).\]
\end{proof}

\begin{lem}\label{lem:pair'det}
a) The limit
\[
\left\langle\det(w),e_{ij}(z)\right\rangle' :=
\lim_{(z_1,w_1)\to (z,w)}
 \frac{\left\langle\det(w_1),e_{ij}(z_1)\right\rangle}
   {\theta(w_1/z_1)\theta(q^2w_1/z_1)\cdots\theta(q^{2(n-2)}w_1/z_1)}
\]
exists for any $z,w\in\C^\times$.

b) We have
\begin{equation}\label{eq:detpair'}
\mu_l\left(\left\langle
\det(w),e_{11}(z)\right\rangle' 1\right)
\det(w)e_{11}(z)=
\mu_r\left(\left\langle
\det(w),e_{11}(z)\right\rangle' 1\right)
e_{11}(z)\det(w)
\end{equation}
for any $z,w\in\C^\times$.
\end{lem}
\begin{proof}
a) We must show that
$\left\langle\det(w),e_{ij}(z)\right\rangle$ vanishes
for $q^{2k}w/z\in p^\Z$, where $k\in\{0,1,\ldots,n-2\}$.
Applying the Laplace expansion
\eqref{eq:rightLaplace} twice we get
\begin{equation}\label{eq:detexpansion}
\det(w)=\sum_{I_1\cup I_2\cup I_3=[1,n]}
 S_r(I_1,I_2,I_3;\la)
 \xi_{I_1}^{J_1}(w)
 \xi_{I_2}^{J_2}(q^{2\#J_1}w)
 \xi_{I_3}^{J_3}(q^{2(\#J_1+2)}w)
\end{equation}
where $J_1=\{1,\ldots,k\}$, $J_2=\{k+1,k+2\}$, $J_3=\{k+3,\ldots,n\}$.
Substituting this in the pairing and applying the multiplication-comultiplication
relation \eqref{eq:pairing_ab_c} we see that each term contains
a factor of the form
$\left\langle\xi_X^Y(q^{2k}w),e_{xy}(z)\right\rangle$ where
$\#X=\#Y=2$ which indeed vanishes for $q^{2k}w/z\in p^\Z$
by Lemma \ref{lem:pairingzeros}.

b) By \eqref{eq:cobraidinghomogen},
$\left\langle\det(w),e_{xy}(z)\right\rangle=0$ if $x\neq y$.
Thus \eqref{eq:detpair'} can be written
\[
\sum_x \mu_l\left(\left\langle\det(w),e_{1x}(z)\right\rangle '1\right)\det(w)e_{x1}(z)
=\sum_x \mu_r\left(\left\langle\det(w),e_{x1}(z)\right\rangle '1\right)e_{1x}(z)\det(w).
\]
If $q^{2k}w/z\notin p^\Z$ for any $k\in\{0,1,\ldots,n-2\}$,
this follows from the cobraiding identity \eqref{eq:cobraidingidentity}
with $a=\det(w), b=e_{11}(z)$ by dividing by the
nonzero number $\prod_{k=0}^{n-2}\theta(q^kw/z)$.

So assume $q^{2k}w/z\in p^\Z$ for some $k\in\{0,1,\ldots,n-2\}$.
We again use the iterated Laplace expansion \eqref{eq:detexpansion}.
For simplicity of notation, we write it
as $\det(w)=\sum a_1a_2a_3$ where $a_2$ is the $2\times 2$ minor.
Put $b=e_{11}(z)$.
Substituting this, and expanding 
$\langle a_1a_2a_3,b'\rangle$ using \eqref{eq:pairing_ab_c},
we get after simplification
\begin{multline*}
\sum_x \mu_l\!\left(\left\langle\det(w),e_{1x}(z)\right\rangle '1\right)\det(w)e_{x1}(z)
=\\
=\frac{1}{\prod_{\substack{m=0\\m\neq k}}^{n-2}\theta(q^{2m}w/z)}\sum_{(a),(b)}
 \mu_l\!\left(\left\langle a_1',b'\right\rangle 1\right) a_1''
 \mu_l\!\left(\big\langle a_2',b''\big\rangle' 1\right)
  a_2''
 \mu_l\!\left(\left\langle a_3',b^{(3)}\right\rangle 1\right) a_3'' b^{(4)} 
\end{multline*}
Now using the cobraiding identity \eqref{eq:cobraidingidentity} and
its primed version for quadratic minors \eqref{eq:keyidentity},
we can move the $b$ all the way to the left. Doing the steps backwards
the claim follows.
\end{proof}

It remains to calculate $\left\langle \det(w),e_{11}(z)\right\rangle' 1$.

\begin{lem} \label{lem:detcalc} We have
\begin{equation}
\left\langle \det(w),e_{11}(z)\right\rangle' 1 = q^{n-1}\th(q^{2n}w/z).
\end{equation}
\end{lem}
\begin{proof}
Expanding $\det(w)$ using
the left expansion formula \eqref{eq:leftminor} with $\si=
\begin{pmatrix}
1 & 2 & \cdots & n \\
n & n-1 & \cdots & 1
\end{pmatrix}$, the longest element in $S_n$,
and applying \eqref{eq:pairing_ab_c} repeatedly we have (putting $I=[1,n]$)
\begin{align*}
\left\langle \det(w), e_{11}(z)\right\rangle &=
\sum_{\tau\in S_n} \left\langle
 \frac{\tau(F_I(\rho))}{\si(F_I(\la))}
  e_{1\tau(n)}(q^{2(n-1)}w)\cdots e_{n\tau(1)}(w), e_{11}(z)
 \right\rangle=\\
&=\sum_{\substack{\tau\in S_n\\x_1,\ldots,x_{n-1}}}
\tau(F_I(\ze))
 \left\langle e_{1\tau(n)}(q^{2(n-1)}w), e_{1x_1}(z) \right\rangle T_{\om(x_1)}
 \cdots \\&\qquad\qquad
 \cdots T_{\om(x_{n-1})}
 \left\langle e_{n\tau(1)}(w), e_{x_{n-1}1}(z)\right\rangle
 \si(F_I(\ze)^{-1}).
\end{align*}

One proves inductively that in all nonzero terms
we have $\tau(j)=n+1-j$ and $x_{n-j}=1$ for all $1\le j \le n-1$
by looking from right to left: 
$\left\langle e_{n\tau(1)}(w),e_{x_{n-1}1}(z)\right\rangle=
\tilde R_{n x_{n-1}}^{\tau(1) 1}(\ze,w/z)T_{-\om(n)-\om(x_{n-1})}$
which, if $1\neq n$, is nonzero only for $\tau(1)=n$ and $x_{n-1}=1$
by \eqref{eq:Rmatrixelts}.
Then looking at the second pairing from the right we see that
$\tau(2)=n-1$ and $x_{n-2}=1$ if it is nonzero, and so on.
Thus only the term $\tau=\si$ and $x_1=\cdots=x_{n-1}=1$
survives and it equals
\begin{multline*}
\si(F_I(\ze))\tilde R_{11}^{11}(\ze,\frac{q^{2(n-1)}w}{z})T_{-\om(1)}
 \tilde R_{21}^{21}(\ze,\frac{q^{2(n-2)}w}{z})T_{-\om(2)}\cdots \\
 \qquad\qquad\cdots
 T_{-\om(n-1)}\tilde R_{n1}^{n1}(\ze,\frac{w}{z})T_{-\om(n)-\om(1)}
 \si(F_I(\ze))^{-1}
\end{multline*}
Using that $\si(F_I(\ze))=\prod_{i<j}E(\ze_{ji}+1)$ and that
$\tilde R_{j1}^{j1}(\ze-\om(1),z)=\tilde\alpha(\ze_{j1}+1,z)=
q\th(z)\frac{E(\ze_{j1}+2)}{E(\ze_{j1}+1}$
we get
\begin{multline*}
 q^{n-1}\theta(q^{2n}w/z)\cdot\theta(q^{2(n-2)}w/z)\theta(q^{2(n-3)}w/z)
 \cdots \theta(w/z)\cdot \\
 \cdot
\prod_{i<j}E(\ze_{ji}+1) 
\prod_{1<j}\frac{E(\ze_{j1}+2)}{E(\ze_{j1}+1)}
T_{-\om(1)}\prod_{i<j}E(\ze_{ji}+1)^{-1}.
\end{multline*}
The factors involving the dynamical variable $\ze$ cancel and the claim
follows.
\end{proof}

By Lemma \ref{lem:pair'det} b) and Lemma \ref{lem:detcalc} we conclude
that $\det(w)$ commutes with $e_{11}(z)$ if $q^{2n}w/z\notin p^\Z$.
By applying an automorphism from the $S_n\times S_n$-action on $A$
as defined in Section \ref{sec:Snaction} and using that 
$\det(z)$ is fixed by those, by relation \eqref{eq:Snonminors},
we conclude that
$\det(w)$ commutes with any $e_{ij}(z)$ as long as $q^{2n}w/z\notin p^\Z$.

For the remaining case we can note that relations
\eqref{eq:rel1-4},\eqref{eq:rel5}
imply that there is a $\C$-linear map $T:\Fell(M(n))\to\Fell(M(n))$
such that $T(ab)=T(b)T(a)$ for all $a,b\in \Fell(M(n))$,
given by
\[T\big( e_{ij}(z)\big)=e_{ij}(z^{-1}),\quad T\big(f(\la)\big)=f(-\la),
\quad T\big(f(\rh)\big)=f(-\rh),\]
for all $f\in M_{\hf^\ast}$, $i,j\in [1,n]$ and $z\in\C^\times$.
One verifies that $T(\det(z))=\det(q^{-2(n-1)}z^{-1})$.

We have proved that $[\det(w),e_{ij}(z)]=0$ if $q^{2n}w/z\notin p^\Z$.
Assume $q^{2n}w/z\in p^\Z$. Then
\[T\big([\det(w),e_{ij}(z)]\big)=[e_{ij}(z^{-1}),\det(q^{-2(n-1)}w^{-1})]=0\]
since $q^{-2(n-1)}w^{-1}/z^{-1}=q^2 (q^{2n}w/z)^{-1}\notin p^\Z$.
This finishes the proof of Theorem \ref{thm:det} b).
\end{proof}


\subsection{The antipode} \label{sec:antipode}
We use the following definition for the antipode, given in \cite{KoRo}.
\begin{dfn}
 An \emph{$\hf$-Hopf algebroid} is an $\hf$-bialgebroid $A$ equipped with a $\C$-linear
 map $S:A\to A$, called the antipode, such that
\begin{equation}
 S(\mu_r(f)a)=S(a)\mu_l(f),\quad S(a\mu_l(f))=\mu_r(f)S(a),\quad a\in A,f\in M_{\hf^*},
\end{equation}
\begin{equation}\label{eq:antipodeaxiom}
\begin{split}
 m\circ (\Id\otimes S)\circ \De(a)&=\mu_l(\ep(a)1),\quad a\in A,\\
 m\circ (S\otimes\Id)\circ \De(a)&=\mu_r(T_\al(\ep(a)1)),\quad a\in A_{\al\be},
\end{split}
\end{equation}
where $m$ denotes the multiplication and $\ep(a)1$ is the result of applying
the difference operator $\ep(a)$ to the constant function $1\in M_{\hf^*}$.
\end{dfn}

Let $\Fell(M(n))[\det(z)^{-1}:z\in\C^\times]$ be the polynomial
algebra in uncountably many variables
$\det(z)^{-1}$, $z\in\C^\times$, with coefficients in
$\Fell(M(n))$.
We define $\Fell(GL(n))$ to be
\[ \Fell(M(n))[\det(z)^{-1}:z\in\C^\times ] / J\]
where $J$ is the ideal generated by the relations
$\det(z)\det(z)^{-1}=1=\det(z)^{-1}\det(z)$ for each $z\in \C^\times$.
We extend the bigrading of $\Fell(M(n))$ to
$\Fell(M(n))[\det(z)^{-1}:z\in\C^\times]$
by requiring that $\det(z)^{-1}$ has bidegree $0,0$
for each $z\in\C^\times$. Then $J$ is homogenous
and the bigrading descends to $\Fell(GL(n))$.
We extend the comultiplication and counit by requiring that $\det(z)^{-1}$ is
grouplike for each $z\in \C^\times$, i.e. that
\[\De(\det(z)^{-1})=\det(z)^{-1}\otimes \det(z)^{-1},\qquad
\ep(\det(z)^{-1})=1.\]
Here $1$ denotes the identity operator in $D_{\hf}$.
One verifies that $J$ is a coideal and that $\ep(J)=0$,
which induces operations $\De,\ep$ on $\Fell(GL(n))$.
In this way $\Fell(GL(n))$ becomes an $\hf$-bialgebroid.
This algebra is nontrivial since $\ep(J)=0$ implies that $J$ is a proper ideal.

For $i\in\oneton$ we set $\hat\imath=\{1,\ldots,i-1,i+1,\ldots,n\}$.
\begin{thm}
 $\Fell(GL(n))$ is an $\hf$-Hopf algebroid with antipode $S$ given by
\begin{equation}\label{eq:antipodedef1}
 S(f(\la))=f(\rh),\qquad S(f(\rh))=f(\la),
\end{equation}
 \begin{equation}\label{eq:antipodedef2}
 S(e_{ij}(z))=\frac{S_r(\hat\jmath,\{j\};\la)}{S_r(\hat\imath,\{i\};\rh)}
 \det(q^{-2(n-1)}z)^{-1}\xi_{\hat\jmath}^{\hat\imath}(q^{-2(n-1)}z),
 \end{equation}
\begin{equation}\label{eq:antipodedef3}
S(\det(z)^{-1})=\det(z),
\end{equation}
for all $f\in M_{\hf^*}$, $i,j\in [1,n]$ and $z\in\C^\times$.
\end{thm}
\begin{proof} We proceed in steps.

\vspace{0.2cm}
\noindent\textbf{Step 1.}\; 
Define $S$ on the generators of $\Fell(M(n))$ by \eqref{eq:antipodedef1},
\eqref{eq:antipodedef2}.
We show that the antipode axiom \eqref{eq:antipodeaxiom} holds if $a$
is a generator. Indeed for $a=f(\la)$ or $a=f(\rh)$, $f\in M_{\hf^*}$ this is
easily checked. Let $a=e_{ij}(z)$. Using the right Laplace expansion \eqref{eq:rightLaplace}
with $J_1=\hat\imath$, $J_2=\{j\}$, $I=\oneton$ and $z$ replaced by $q^{-2(n-1)}z$ we obtain 
\begin{equation}\label{eq:Sproof1}
\sum_{x=1}^n S(e_{ix}(z))e_{xj}(z)=\de_{ij}.
\end{equation}
Similarly, using the left Laplace expansion \eqref{eq:leftLaplace} with
$I_1=\{i\}$, $I_2=\hat\jmath$, $J=\oneton$ and $z$ replaced by $q^{-2(n-1)}z$,
together with the identity \eqref{eq:SlSrrelation}, we get
\begin{equation}\label{eq:Sproof2}
 \sum_{x=1}^n e_{ix}(z)S(e_{xj}(z))=\de_{ij},
\end{equation}
using also the crucial fact that, by Theorem \ref{thm:det}, 
$e_{ij}(z)$ commutes in $\Fell(M(n))$ with
$\det(q^{-2(n-1)}z)$ and hence in $\Fell(GL(n))$ with $\det(q^{-2(n-1)}z)^{-1}$.
This proves that the antipode axiom \eqref{eq:antipodeaxiom} is satisfied for $a=e_{ij}(z)$.

\vspace{0.2cm}
\noindent\textbf{Step 2.}\; 
We show that $S$ extends to
a $\C$-linear map $S:\Fell(M(n))\to \Fell(GL(n))$ satisfying
$S(ab)=S(b)S(a)$. For this we must verify that $S$ preserves the relations,
\eqref{eq:FellMnscalarrel},\eqref{eq:RLLmatrixelements}, \eqref{eq:rel5}
of $\Fell(M(n))$.
Since $S(e_{ij}(z))\in\Fell(GL(n))_{\om(\hat\jmath),\om(\hat\imath)}$
and $\om(i)+\om(\hat\imath)=0$, we have
\begin{multline*}
S(e_{ij}(z))S(f(\la))=S(e_{ij}(z))f(\rh)=f(\rh-\om(\hat\imath))S(e_{ij}(z))=\\=
f(\rh+\om(i))S(e_{ij}(z))=S\big(f(\la+\om(i))\big)S(e_{ij}(z))
\end{multline*}
similarly, $S(e_{ij}(z))S(f(\rh))=S\big(f(\rh+\om(j))\big)S(e_{ij}(z))$ so relations
\eqref{eq:FellMnscalarrel} are preserved.
Next, consider the RLL relation
\begin{equation}\label{eq:Sproof3}
\sum_{x,y=1}^n R_{ac}^{xy}(\la,\frac{z_1}{z_2})e_{xb}(z_1)e_{yd}(z_2)=
\sum_{x,y=1}^n R_{xy}^{bd}(\rh,\frac{z_1}{z_2})e_{cy}(z_2)e_{ax}(z_1).
\end{equation}
Multiply \eqref{eq:Sproof3} from the left by $S(e_{ic}(z_2))$ and from the right by $S(e_{dk}(z_2))$,
sum over $c,d$ and use \eqref{eq:Sproof1},\eqref{eq:Sproof2} to obtain
\[\sum_{x,c}R_{ac}^{xk}(\la-\om(\hat c),\frac{z_1}{z_2})S(e_{ic}(z_2))e_{xb}(z_1)=
\sum_{x,d}R_{xi}^{bd}(\rh-\om(\hat\imath),\frac{z_1}{z_2})e_{ax}(z_1)S(e_{dk}(z_2)).\]
Then multiply from the left by $S(e_{ja}(z_1))$ and from the right by $S(e_{bl}(z_1))$,
sum over $a,b$ and use \eqref{eq:Sproof1},\eqref{eq:Sproof2} again to get
\begin{multline}\label{eq:Sproof4}
 \sum_{a,c}R_{ac}^{lk}(\la-\om(\hat a)-\om(\hat c),\frac{z_1}{z_2})S(e_{ja}(z_1))S(e_{ic}(z_2))=\\
 =\sum_{b,d}R_{ji}^{bd}(\rh-\om(\hat\jmath)-\om(\hat\imath),\frac{z_1}{z_2})S(e_{dk}(z_2))S(e_{bl}(z_1)).
\end{multline}
Since $S(e_{ij}(z))\in\Fell(GL(n))_{\hat\jmath,\hat\imath}$ and
$R_{ji}^{bd}(\rh-\om(\hat\jmath)-\om(\hat\imath),\frac{z_1}{z_2})=
R_{ji}^{bd}(\rh-\om(\hat b)-\om(\hat d),\frac{z_1}{z_2})$ by the $\hf$-invariance of $R$,
\eqref{eq:Sproof4} can be rewritten
\[ \sum_{a,c}S(e_{ja}(z_1))S(e_{ic}(z_2))R_{ac}^{lk}(\la,\frac{z_1}{z_2})
 =\sum_{b,d}S(e_{dk}(z_2))S(e_{bl}(z_1))R_{ji}^{bd}(\rh,\frac{z_1}{z_2}).\]
This is the result of formally applying $S$ to the RLL-relations, proving that
$S$ preserves \eqref{eq:RLLmatrixelements}. Similarly \eqref{eq:rel5} is preserved.

\vspace{0.2cm}
\noindent\textbf{Step 3.}\; 
Since, by the above steps, \eqref{eq:antipodeaxiom} holds on the generators of $\Fell(M(n))$
and $S(ab)=S(b)S(a)$ for all $a,b\in \Fell(M(n))$, it follows
that \eqref{eq:antipodeaxiom} holds for any $a\in \Fell(M(n))$.
By taking in particular $a=\det(z)$ we get
\[\det(z)S(\det(z))=1\qquad\text{and}\qquad S(\det(z))\det(z)=1\]
respectively. Thus, definining $S$ on $\det(z)^{-1}$ by \eqref{eq:antipodedef3},
the relations $\det(z)\det(z)^{-1}=1=\det(z)^{-1}\det(z)$ are preserved by $S$.
Hence $S$ extends to an anti-multiplicative $\C$-linear map
$S:\Fell(GL(n))\to\Fell(GL(n))$ satisfying the antipode axiom \eqref{eq:antipodeaxiom}
on $\Fell(M(n))$ and on $\det(z)^{-1}$. Hence \eqref{eq:antipodeaxiom}
holds for any $a\in\Fell(GL(n))$.
\end{proof}

\section{Concluding remarks} \label{sec:discussion}

To define the antipode we only needed that $e_{ij}(z)$ commutes
with $\det(q^{-2(n-1)}z)$. This can also be proved using the Laplace expansions.

Perhaps one could avoid problems with spectral poles and zeros of the 
R-matrix by thinking of the algebra as generated by meromorphic sections
of a $M_{\hf^\ast\oplus \hf^\ast}$-line bundle over the elliptic curve $\C^\times / \{z\sim pz\}$.
In this direction we found that the relation $e_{ij}(pz)=q^{\la_i-\rh_j}e_{ij}(z)$
respects the $RLL$-relation (here $\hf$ should be the Cartan subalgebra of 
$\mathfrak{gl}_n$). This relation should then most likely be added to
the algebra.

It would be interesting to
develop harmonic analysis for the elliptic $GL(n)$ quantum group,
similarly to \cite{KoRo}. 
In this context it is valuable to have an abstract algebra to work with,
and not only a tensor category analogous to a category of representations.
For example the analog of the Haar measure seems most naturally defined
as a certain linear functional on the algebra.

\section*{Acknowledgements}
The author is greatly indebted to H. Rosengren
for many inspiring and helpful discussions, and to J. Stokman for his
support and helpful comments. 
The author is also grateful to an anonymous referee
of an earlier version of this paper
for detailed comments which has led to improvement of the section
on centrality of the determinant.

\noindent\textsc{
University of Amsterdam,
Korteweg-de Vries Institute for Mathematics,
Postbus 94248,
1090 GE Amsterdam,
Netherlands}

\noindent\texttt{Email: jonas.hartwig@gmail.com}

\noindent\texttt{URL: http://staff.science.uva.nl/{\textasciitilde}jhartwig}

\end{document}